\documentclass[a4paper,10pt]{amsart}
\usepackage{amssymb}
\usepackage{latexsym}
\usepackage{amsmath}
\usepackage{enumerate}

\usepackage{amsmath, hyperref}
\newtheorem{theorem}{Theorem}[section]

\theoremstyle{definition}
\newtheorem{definition}[theorem]{Definition}

\newcommand{\PA}{\mathbf{PA}}

\numberwithin{equation}{section}

\title{On the depth of G\"{o}del's incompleteness theorem}
\author{Yong Cheng}
\address{School of Philosophy, Wuhan University, China}
\email{world-cyr@hotmail.com}
\begin{document}

\begin{abstract}
In this paper, we use G\"{o}del's incompleteness theorem as a case
study for investigating  mathematical depth.
We take for granted the widespread judgment by mathematical logicians that G\"{o}del's incompleteness theorem is deep, and focus on the philosophical question of
what its depth consists in.  We focus on the methodological study of the depth of G\"{o}del's incompleteness theorem, and propose three criteria to account for its depth:  influence, fruitfulness, and unity. Finally, we give some explanations for our account of the depth of G\"{o}del's incompleteness theorem.
\end{abstract}

\maketitle

\section{Introduction}

``Mathematical depth" is an often used notion  when mathematicians assess and
evaluate the work of their fellows. Mathematics is full of disagreements over what is deep work.  And even mathematicians who do agree in judging one thing or another deep, are likely to disagree sharply on what makes it deep.
The notion of mathematical depth is not   well-defined, and there is no single widely accepted account of mathematical depth.
There is a lot of discussion of mathematical depth from the literature (see [Ernst, Heis, Maddy, McNulty and Weatherall, 2015a, 2015b][Arana, 2015] [Gray, 2015] [Lange, 2015][Stillwell, 2015][Urquhart, 2015]).
The following four types of questions  have been widely discussed in the current philosophical investigation (see [Ernst, Heis, Maddy, McNulty and Weatherall, 2015a, pp.155-162]):
(1) is there agreement that the cited examples are deep or not deep? (2) Are there commonalities in the kinds of features cited in defense of depth and
non-depth assessments in the various examples? (3) Is depth the same as or different from such notions as fruitfulness,
surprisingness, importance, elegance, difficulty, fundamentalness,
explanatoriness, beauty, etc.? (4)  Is depth an objective feature or something essentially tied to our interests,
abilities, and so on?

%four accounts of mathematical depth will be considered.

In this paper, the bearer of depth is mathematical theorems.
There may be different ways of being deep, and different theorems may have distinct criteria for  their depth.
It is hard for us to give a universal account of mathematical depth.  In this paper, we do not attempt to give a universal  account of mathematical depth. Instead, our strategy in this paper is to use G\"{o}del's incompleteness theorem as a case
study for investigating mathematical depth.
In this paper, we take for granted the widespread judgment by  mathematical logicians  that G\"{o}del's incompleteness theorem is deep, and focus on the philosophical question of
what its depth consists in.  In particular, we focus on the methodological study of the depth of G\"{o}del's incompleteness theorem: we attempt to find some fundamental criteria to account for the depth of the incompleteness theorem based on the current research on incompleteness  from the literature.
After briefly introducing G\"{o}del's incompleteness theorem,  we account for the depth of the incompleteness theorem according to the following three criteria: influence, fruitfulness, and unity.
In this paper, we  only cover  the most important mathematical evidences (as far as we know) of  the three criteria. We make no attempt at completeness here. There are many more mathematical evidences of the depth of G\"{o}del's incompleteness theorem from the  literature than we can cover here. This paper is a modest attempt to bring some coherence to philosophical understandings of the depth of G\"{o}del's incompleteness theorem.
This work may help us  to get a sense of  the viability of a methodological inquiry into what
counts as deep mathematics and why.
We hope this work will advance our philosophical understanding of mathematical depth, perhaps making way for more unified accounts of mathematical depth to follow.

%Instead, the motivation of this paper is to  account for the depth of a specific theorem in mathematical logic which  has been judged deep by many  mathematical logicians: G\"{o}del's incompleteness theorem.
%
%

This paper is written in the same spirit of Arana's [2015]. After introducing  Szemer\'{e}di's Theorem and gesturing at its proofs, Arana [2015] articulates four different accounts of the depth of Szemer\'{e}di's Theorem (genetic views, evidentialist views, consequentialist views, and cosmological
views), and  indicates ways in which each is apt and inapt for characterizing depth as it occurs
in mathematical practice.
The genetic view of depth  identifies a deep theorem as one
proved by sufficiently talented mathematicians.
The evidentialist view of depth  links the depth of a theorem with some quality of its proof.
The consequentialist view of depth
measures the depth of a theorem by some quality of its consequences, or of the consequences of its proofs.
The cosmological view of depth measures the depth of a theorem by measuring the order the theorem has established and the unexpected
structure it has revealed.
However, Arana's four views of depth have features
that many philosophers would consider problematic: either by falling into vagueness, by failing to pick out theorems that obviously count as deep, or by making depth subjective (see [Ernst, Heis, Maddy, McNulty, and Weatherall, 2015a]).
Our three criteria of the depth of G\"{o}del's incompleteness theorem are isolated from Arana's  evidentialist, consequentialist, and cosmological views of depth. We basically agree with Arana's analysis of these views of depth.
The focus of this paper is the justification of the depth of G\"{o}del's incompleteness theorem from the three criteria we propose. %In our account of the depth of G\"{o}del's incompleteness theorem, we propose specific mathematical evidences of these three criteria which indicates

%Arana's four views of depth focus on what aspects of a theorem are linked with its depth. In this paper, we do not intend to discuss the general criteria of mathematical depth; instead we propose three specific criteria to justify the depth of G\"{o}del's incompleteness theorem.
%Compare our three criteria with Arana's four views of depth, our three criteria are more specific than Arana's four views of depth.

In this paper, $\sf G1$ stands for G\"{o}del's first incompleteness theorem, and $\sf G2$ stands for G\"{o}del's second incompleteness theorem.
This paper is structured as follows. In Section 1, we introduce the motivation and main content of this paper. In Section 2, we give a brief overview of G\"{o}del's incompleteness theorem and a sketch of main ideas of the  proof of $\sf G1$ and $\sf G2$.
In Section 3,  we account for  the depth of G\"{o}del's incompleteness theorem based on the following three criteria: influence, fruitfulness and unity.
In Section 4, we give some explanations for our account of the depth of G\"{o}del's incompleteness theorem.

%based on known research on mathematical depth  from the literature
%The criterion measures laboriousness in terms of several properties of proofs: length, complexity, difficulty, and %impurity.

%Let us consider another evidentialist strategy: a theorem is deep if it has many different proofs.

%From \cite{Stillwell}, A deep result emerges or comes to light only after the
%intermediate layers of concepts and proofs have been elucidated;

%\section{Preliminaries}\label{pre sec}

\section{G\"{o}del's incompleteness theorem}

In this section, we give an overview of G\"{o}del's incompleteness theorem. For textbooks on G\"{o}del's incompleteness theorem, we refer to [Enderton 2001], [Murawski, 1999], [Lindstr\"{o}m, 1997], [Smith, 2007],  [Boolos, 1993].
For survey papers on G\"{o}del's incompleteness theorem, we refer to [Smory\'{n}ski, 1977], [Beklemishev, 2010], [Kotlarski, 2004], [Visser, 2016], [Cheng, 2019c].

We first review some basic notions used in this paper.
We focus on  first order theories with a countable language. For a given theory $T$, let \emph{$L(T)$} denote  the language of $T$.
In this paper, arithmetization refers to the method in mathematical logic that replaces  reasonings on the expressions of first order  language by reasonings on natural numbers. For this purpose, the replacement is constructed by some one-to-one mapping of the set of all expressions (in the alphabet of the language under consideration) into the natural number sequence. Relations and operations defined on expressions are transformed by this mapping into relations and operations on natural numbers.  If unless stated otherwise, we always assume the arithmetization of the base theory with a recursive set of non-logical constants. For more details about arithmetization, we refer to [Murawski, 1999].
Under arithmetization, any formula (or finite sequence of formulas)  can be coded by a natural number (called \emph{G\"{o}del's number}). We use \emph{$\ulcorner\phi\urcorner$} to denote the corresponding numeral
of the G\"{o}del number of a formula $\phi$. Given a theory $T$, we say a sentence $\phi$ in $L(T)$ is independent of $T$ if $T\nvdash \phi$ and $T\nvdash \neg\phi$.
A theory $T$ is
\emph{incomplete} if there is a sentence $\phi$ in $L(T)$ such that $\phi$  is independent of $T$; otherwise, $T$ is
\emph{complete}.  A theory $T$ is
\emph{recursively axiomatizable} if it has a recursive set of axioms, i.e.~ the set of G\"{o}del numbers of axioms of $T$ is recursive.
A $n$-ary relation $R(x_1, \cdots, x_n)$ on $\mathbb{N}^n$ is \emph{representable} in a theory $T$ if there is a formula $\phi(x_1, \cdots, x_n)$ such that if $R(m_1, \cdots, m_n)$ holds, then $T\vdash \phi(\overline{m_1}, \cdots, \overline{m_n})$;\footnote{For  $n\in \mathbb{N}$, $\bar{n}$ denotes the corresponding numeral in $L(T)$ for $n$.} and if $R(m_1, \cdots, m_n)$ does not hold, then $T\vdash \neg\phi(\overline{m_1}, \cdots, \overline{m_n})$.\footnote{In this case, we say the formula $\phi(x_1, \cdots, x_n)$ represents the relation $R(x_1, \cdots, x_n)$.}

Robinson Arithmetic $\mathbf{Q}$ is introduced in [Tarski, Mostowski and Robinson, 1953]  as a base axiomatic theory for investigating incompleteness and undecidability.

\begin{definition}\label{def of Q}
Robinson Arithmetic $\mathbf{Q}$  is  defined in   the language $\{\mathbf{0}, \mathbf{S}, +, \times\}$ with the following axioms:
\begin{description}
  \item[$\mathbf{Q}_1$] $\forall x \forall y(\mathbf{S}x=\mathbf{S} y\rightarrow x=y)$;
  \item[$\mathbf{Q}_2$] $\forall x(\mathbf{S} x\neq \mathbf{0})$;
  \item[$\mathbf{Q}_3$] $\forall x(x\neq \mathbf{0}\rightarrow \exists y (x=\mathbf{S} y))$;
  \item[$\mathbf{Q}_4$]  $\forall x\forall y(x+ \mathbf{0}=x)$;
  \item[$\mathbf{Q}_5$] $\forall x\forall y(x+ \mathbf{S} y=\mathbf{S} (x+y))$;
  \item[$\mathbf{Q}_6$] $\forall x(x\times \mathbf{0}=\mathbf{0})$;
  \item[$\mathbf{Q}_7$] $\forall x\forall y(x\times \mathbf{S} y=x\times y +x)$.
\end{description}
\end{definition}

The theory $\mathbf{PA}$ consists of axioms $\mathbf{Q}_1$-$\mathbf{Q}_2$, $\mathbf{Q}_4$-$\mathbf{Q}_7$ in Definition \ref{def of Q} and the axiom scheme of induction $(\phi(\mathbf{0})\wedge \forall x(\phi(x)\rightarrow \phi(\mathbf{S} x)))\rightarrow \forall x  \phi(x)$, where $\phi$ is a formula with at least one free variable $x$.
Let $\mathfrak{N}=\langle\mathbb{N}, +, \times\rangle$ denote the standard model of arithmetic. We say $\phi\in L(\mathbf{PA})$ is a true sentence of arithmetic if $\mathfrak{N}\models\phi$.

We introduce a hierarchy of $L(\mathbf{PA})$-formulas called the \emph{arithmetical hierarchy} (see [P. H\'{a}jek  and P. Pudl\'{a}k, 1993]). \emph{Bounded formulas}   ($\Sigma^0_0$, or $\Pi^0_0$, or $\Delta^0_0$ formula) are built from atomic
 formulas using only propositional connectives and bounded quantifiers (in the form $\forall x\leq y$ or $\exists x\leq y$). A formula is  $\Sigma^0_{n+1}$ if it has the
form $\exists x\phi$ where $\phi$ is $\Pi^0_{n}$. A formula is $\Pi^0_{n+1}$ if it has the form $\forall x\phi$ where $\phi$ is $\Sigma^0_{n}$. Thus, a $\Sigma^0_{n}$-formula has a block of $n$ alternating quantifiers, the first one
being existential, and this block is followed by a bounded formula. Similarly
for $\Pi^0_{n}$-formulas. A formula is $\Delta^0_n$ if it is equivalent to both a $\Sigma^0_{n}$ formula and a $\Pi^0_{n}$ formula in $\mathbf{PA}$.

A theory $T$ is said to be $\omega$-consistent if there is no formula $\varphi(x)$ such that $T \vdash\exists x \varphi(x)$ and for any $n\in\mathbb{N}$, $T \vdash\neg\varphi(\bar{n})$; $T$ is $1$-consistent if there is no such a $\Delta^0_1$ formula $\varphi(x)$.
We say a theory $T$ is $\Sigma^0_{n}$-definable if there is a  $\Sigma^0_{n}$ formula $\alpha(x)$ such that $n$ is the G\"{o}del number of some sentence of $T$ if and only if $\mathfrak{N}\models \alpha(\overline{n})$. A theory $T$ is $\Sigma^0_{n}$-sound if  for all $\Sigma^0_{n}$ sentences $\phi$,  $T\vdash\phi$ implies $\mathfrak{N}\models\phi$.
A theory $T$ is $\Pi^0_{n}$-decisive if  for all $\Pi^0_{n}$ sentences $\phi$, either $T\vdash\phi$ or $T \vdash\neg\phi$ holds.

The notion of interpretation provides us with a method to measure and compare the strength of different theories in different languages.
Informally, an interpretation of a theory $T$ in a theory $S$ is a mapping from formulas
of $T$ to formulas of $S$ that maps all axioms of $T$ to
sentences provable in $S$.
For the precise
definition of interpretation, we refer to [Visser, 2011] for more details.
Let $T\unlhd S$ denote that $T$ is interpretable in $S$, and let $T\lhd S$ denote that $T\unlhd S$ but $S\unlhd T$ does not hold. In this paper, we say that $T$ is weaker  than $S$ w.r.t.~ interpretation if $T\lhd S$.

G\"{o}del proves his incompleteness theorem in [G\"{o}del, 1931] for a certain formal
system $\mathbf{P}$ related to Russell-Whitehead's Principia Mathematica and based on the
simple theory of types over the natural number series and the Dedekind-Peano
axioms (see [Beklemishev, 2010, p.3]).
G\"{o}del's original first incompleteness theorem ([G\"{o}del, 1931]) says that for
any formal theory $T$ formulated in the language of $\mathbf{P}$ and obtained by adding a primitive recursive set of axioms to the system $\mathbf{P}$, if $T$ is $\omega$-consistent, then $T$ is incomplete.
The following  is a modern reformulation of G\"{o}del's first incompleteness theorem.

\begin{theorem}[G\"{o}del's first incompleteness theorem $(\sf G1)$]\label{G1}~
Let $T$ be a recursively axiomatized extension of $\mathbf{PA}$. Then there exists a  G\"{o}del's sentence $\mathbf{G}$ such that:
\begin{itemize}
   \item if $T$ is consistent, then $T\nvdash \mathbf{G}$;
  \item if $T$ is $\omega$-consistent, then $T\nvdash \neg\mathbf{G}$.
\end{itemize}
\end{theorem}

From Theorem \ref{G1},  if $T$ is $\omega$-consistent, then G\"{o}del's sentence $\mathbf{G}$ is independent of $T$ and hence
$T$ is incomplete.

Now, we give a sketch of the main idea of G\"{o}del's proof of $\sf G1$.
In the rest of this section, we assume that $T$ is a recursively axiomatized consistent extension of $\mathbf{PA}$ in $L(\mathbf{PA})$.
G\"{o}del's proof of  the incompleteness theorem depends on a long chain of ideas, each involving a significant insight.
The three main ideas in G\"{o}del's proof of $\sf G1$ are arithmetization of the syntax of $T$,  representability of recursive functions in $\mathbf{PA}$ and  self-reference construction of G\"{o}del's sentence.\footnote{Before G\"{o}del, Emil Post independently discovered a statement undecidable within Principia
Mathematica whose truth could nevertheless be established by
metamathematical considerations. But  he never published the result. See [Murawski, 1999, p. 203].}
Under the arithmetization, we could establish the one-to-one correspondence between expressions of $L(T)$ and natural numbers. Under this correspondence, we can translate metamathematical
statements about the formal theory $T$ into statements
about natural numbers. Moreover, fundamental metamathematical relations can be translated in this way into
certain recursive relations, hence into relations representable in the
theory $T$. Consequently, one can speak about a formal system of arithmetic and about its properties
as a theory in the system itself (see [Murawski, 1999])! This is the essence of G\"{o}del's
idea of arithmetization.\footnote{For more details on arithmetization, we refer to [Murawski, 1999].}

Now, we can define some relations on $\mathbb{N}$ which express some metamathematical properties of $T$.
For example, we can define a binary relation on $\mathbb{N}^2$ as follows: $Proof_T(m,n)$ iff $n$ is the G\"{o}del's number of a proof in $T$ of the formula with G\"{o}del number $m$.
Moreover, we can prove that the  relation $Proof_T(m,n)$  is
recursive. G\"{o}del proves that every recursive relation is representable
in $\PA$. Let  $\mathbf{Proof}_T(x,y)$ be the formula  which represents $Proof_T(m,n)$ in $\mathbf{PA}$.\footnote{One can speak about the property of $T$  in $\mathbf{PA}$ itself via  arithmetization  and representability!} From the formula $\mathbf{Proof}_T(x,y)$, we can define the provability predicate $\mathbf{Prov}_T(x)$ as $\exists y \mathbf{Proof}_T(x,y)$. Finally, G\"{o}del effectively constructes  G\"{o}del's sentence $\mathbf{G}$ which asserts its own unprovability in $T$ (i.e.~ $T \vdash \mathbf{G} \leftrightarrow\neg \mathbf{Prov}_T(\ulcorner \mathbf{G} \urcorner))$.
G\"{o}del shows  that  if $T$ is consistent, then $T\nvdash\mathbf{G}$; and if $T$ is $\omega$-consistent, then $T\nvdash\neg\mathbf{G}$.

Since we will discuss general provability predicates based on proof predicates, now we give a general definition of proof predicate which is a generalization of properties of G\"{o}del's proof predicate $Proof_T(x,y)$.
We say a formula $\mathbf{Prf}_T(x, y)$ is a \emph{proof predicate} of $T$ if it satisfies the following conditions:
\begin{enumerate}[(1)]
  \item $\mathbf{Prf}_T(x, y)$ is $\Delta^0_1$;
  \item $\mathbf{PA} \vdash \forall x(\mathbf{Prov}_T(x) \leftrightarrow\exists y \mathbf{Prf}_T(x, y))$;
  \item  for any $n \in\omega$ and formula $\phi, \mathbb{N}\models \mathbf{Proof}_T(\ulcorner\phi\urcorner, \overline{n}) \leftrightarrow \mathbf{Prf}_T(\ulcorner\phi\urcorner, \overline{n})$;
      \item $\mathbf{PA} \vdash \forall x\forall x^{\prime} \forall y (\mathbf{Prf}_T(x, y) \wedge \mathbf{Prf}_T(x^{\prime}, y) \rightarrow x = x^{\prime})$.
\end{enumerate}

Note that each proof predicate represents the relation ``$y$ is the code of a proof in $T$ of a formula with G\"{o}del number $x$". We  define the  provability predicate $\mathbf{Pr}_T(x)$ from a proof predicate $\mathbf{Prf}_T(x,y)$ by $\exists y\, \mathbf{Prf}_T(x,y)$, and   the  consistency statement $\mathbf{Con}(T)$ from a provability predicate $\mathbf{Pr}_T(x)$ by $\neg \mathbf{Pr}_T (\ulcorner \mathbf{0}\neq\mathbf{0}\urcorner)$.

The following conditions $\mathbf{D1}$-$\mathbf{D3}$ are called drivability conditions of provability predicate $\mathbf{Pr}_T(x)$.
\begin{description}
  \item[D1] If $T \vdash\varphi$, then $T \vdash \mathbf{Pr}_T(\ulcorner\varphi\urcorner)$;
  \item[D2] $T \vdash\mathbf{Pr}_T(\ulcorner\varphi\urcorner)\rightarrow(\mathbf{Pr}_T(\ulcorner\varphi
      \rightarrow\psi\urcorner)\rightarrow \mathbf{Pr}_T(\ulcorner\psi\urcorner))$;
  \item[D3] $T\vdash \mathbf{Pr}_T(\ulcorner\varphi\urcorner) \rightarrow \mathbf{Pr}_T(\ulcorner\mathbf{Pr}_T(\ulcorner\varphi\urcorner)\urcorner)$.
\end{description}
We say a provability predicate is standard if it satisfies conditions $\mathbf{D1}$-$\mathbf{D3}$.
One important non-standard provability predicate is
Rosser provability predicate $\mathbf{Pr}_T^R(x)$ introduced by
Rosser [1936] to improve G\"{o}del's first incompleteness theorem. The Rosser provability predicate $\mathbf{Pr}^R_{T}(x)$ is defined as the formula $\exists y(\mathbf{Prf}_{T}(x,y)\wedge \forall z\leq y\neg \mathbf{Prf}_{T}(\dot{\neg}(x),z))$, where $\dot{\neg}$ is a function symbol expressing a primitive recursive function calculating the code of $\neg \phi$ from the code of $\phi$.

In this paper, unless stated otherwise, we assume that the provability predicate $\mathbf{Pr}_T(x)$ is standard, and $\mathbf{Con}(T)$ defined as $\neg \mathbf{Pr}_T(\ulcorner \mathbf{0}\neq\mathbf{0}\urcorner)$ is the canonical consistence statement of $T$  formulated via a standard provability predicate $\mathbf{Pr}_T(x)$. The importance of standard provability predicate and canonical consistence statement lies in that $\sf G2$ holds for the canonical consistence statement formulated via the standard provability predicate as we will show. However, $\sf G2$ may fail for non-standard provability predicates and non-canonical consistence statements. There are a lot of research on non-standard provability predicates and non-canonical consistence statements from the literature (see [Feferman, 1960] and [Visser, 2011]).

G\"{o}del announces the second incompleteness theorem  in an abstract published in
October 1930: no consistency proof of systems such as Principia, Zermelo-Fraenkel
set theory, or the systems investigated by Ackermann and von Neumann is possible
by methods that can be formulated in these systems (see [Zach, 2007, p.431]). In the modern formulation,
the second incompleteness theorem  states that if $T$ is
consistent, then the arithmetical formula $\mathbf{Con}(T)$ that expresses the consistency of $T$ is not provable in $T$.\footnote{Even if  $\mathbf{Con}(T)$ is not provable in $T$,  $T$ is reflective: for each finite sub-theory $S$ of  $T$, $T \vdash \mathbf{Con}(S)$ (see [Murawski, 1999, Theorem 2.6.12]).}
For the proof of $\sf G2$, we use the key fact that  the provability predicate $\mathbf{Prov}_T(x)$ is standard and satisfies conditions $\mathbf{D1}$-$\mathbf{D3}$. Based on this fact, we can show that $T\vdash \mathbf{Con}(T)\leftrightarrow \mathbf{G}$. Thus, $\sf G2$ holds: if $T$ is consistent, then $T\nvdash \mathbf{Con}(T)$.

For G\"{o}del's proof of $\sf G1$, only assuming that $T$ is consistent  is not enough to show that  G\"{o}del's sentence is independent of $T$. In fact, the optimal condition to show that G\"{o}del's sentence is independent of $T$ is that $T+\mathbf{Con}(T)$ is consistent (see [Isaacson, 2011, Theorems 35-36]).\footnote{This optimal condition is much weaker than $\omega$-consistency.}
Only assuming that $T$ is consistent is not sufficient to show that $T\nvdash\neg\mathbf{Con}(T)$. But we can prove that $\mathbf{Con}(T)$ is independent of $T$ by assuming that $T$ is 1-consistent which is stronger than ``$T$ is consistent".
For more details of proofs of $\sf G1$ and $\sf G2$, we refer to Chapter 2 in [Murawski, 1999].

\section{On the depth of G\"{o}del's incompleteness theorem}

After introducing G\"{o}del's incompleteness theorem, now we turn to the analysis of the depth of G\"{o}del's theorem.
In this section, we account for  the depth of G\"{o}del's incompleteness theorem from the following three criteria: influence,  fruitfulness, and unity.

\subsection{Influence}\label{section on influence}

%Influence on philosophy (mind and machine, G\"{o}del's disjunctive thesis, Hilbert's program), mathematics (concrete incompleteness) and computer science.

In this section, we justify for the influence of G\"{o}del's incompleteness theorem from its impact  on foundations of mathematics, philosophy, mathematics, and theoretic computer science, that is revealed by the  research practice on incompleteness after G\"{o}del.
We make no attempt to exhaustively discuss the full impact of G\"{o}del's incompleteness
theorem and all of the ongoing important research programs that it suggests.

G\"{o}del's incompleteness theorem is one of
the most remarkable and profound discoveries in the 20th century, an important milestone in the history of modern logic.
G\"{o}del's incompleteness theorem has had wide and profound influence on the development of logic,  philosophy, mathematics,  theoretical computer science and other fields, substantially shaping foundations  of
mathematics after 1931.
On the impact of G\"{o}del's incompleteness
theorems, Feferman said: ``their relevance to mathematical
logic (and its offspring in the theory of computation) is paramount; further, their philosophical relevance is significant, but in just what way is far from
settled; and finally, their mathematical relevance
outside of logic is very much unsubstantiated but
is the object of ongoing, tantalizing efforts" (see [Feferman, 2006, p.434]).

The influence of G\"{o}del's incompleteness theorem on  foundations of mathematics is reflected in the following five aspects: (1)  G\"{o}del's incompleteness theorem reveals the independence phenomenon which is common in mathematics and logic;
(2) G\"{o}del's incompleteness theorem shows certain weaknesses and the essential limitation of one given formal system (or the limit of proof and computation);
 (3)  G\"{o}del's incompleteness theorem reveals the essential difference between  the notion of ``provability in $\PA$" and the notion of ``truth in the standard model of arithmetic";\footnote{I.e.~ there is a true sentence of arithmetic which is independent of $\PA$.} (4) G\"{o}del's incompleteness theorem is a blow to Whitehead-Russell's program for proving that all mathematics (or at least quite a lot of it) could be derived solely from logic
in their three-volume
Principia Mathematica; (5) G\"{o}del's incompleteness theorem has profound influence on the development of Hilbert's program.

%Hilbert's program has two main goals: (1) proving all true mathematical propositions in the formalization of mathematics; (2) finding a proof that no contradiction can be obtained in the formalism of mathematics using only finitistic reasoning about definite mathematical objects.

There is extensive literature about the development of Hilbert's program after G\"{o}del's incompleteness theorem, and its effect on mathematical logic (especially proof theory) and philosophy of mathematics (see [Feferman, 1988], [Franks, 2009], [Murawski, 1999], [Simpson, 1988], [Zach, 2007]). The above aspects (2)-(5) are well known, and we only give some explanations of the independence phenomenon  in mathematics and logic.

Nowadays, independence is ubiquitous in logic. The independence phenomenon reveals the big gap between mathematical truth and provability in formal systems. G\"{o}del's incompleteness theorem firstly reveals the independence phenomenon of formal systems at the level of  arithmetic, and shows the essential limitation of any formal system containing ``enough" information of arithmetic. After G\"{o}del, people have found many arithmetic sentences  from classic mathematics that are independent of $\mathbf{PA}$.
The sequent research after G\"{o}del reveals the independence phenomenon of stronger formal systems such as higher order arithmetic and $\mathbf{ZFC}$. For example, G\"{o}del and Cohen show that Continuum Hypothesis $(\mathbf{CH})$ is independent of $\mathbf{ZFC}$, that means $\mathbf{ZFC}$ is essentially incomplete to capture all set-theoretic truth.\footnote{ G\"{o}del proposes the research program to find new axioms of set theory to capture all set-theoretic truth, that is one of the central  motivations of current research in set theory.} Moreover, after G\"{o}del people have found many examples of statements that are independent of $\mathbf{ZFC}$ from varied fields of mathematics such as analysis, algebra,  topology and  mathematical logic.

Now, we give a brief account of the impact of G\"{o}del's incompleteness theorem in mathematics. The incompleteness theorem and its proofs are strikingly original mathematics.
It has been often thought that G\"{o}del's proof of $\sf G1$ uses pure logical method and has no relevance with mathematics: G\"{o}del's sentence constructed via the meta-mathematical method is artificial (not natural), and has no  real mathematical content.
As Isaacson [1987] points out, G\"{o}del's sentence is constructed not by reflecting about arithmetical properties of natural numbers, but by
reflecting about an axiomatic system in which those properties are
formalized. A natural question after G\"{o}del is: can we find examples of natural independent sentences  with  real mathematical contents? We call the incompleteness phenomenon  revealed by  natural independent sentences  with  real mathematical contents as concrete incompleteness.

In fact, the incompleteness phenomenon  prevails in classic mathematics. The research program of concrete incompleteness seeks for natural independent sentences with real mathematical contents.
After G\"{o}del, many natural independent arithmetic sentences with real mathematical contents have been found. These independent sentences have a clear mathematical flavor, and do not refer to the arithmetization of syntax and provability predicate.
Paris and Harrington [1977] proposes  a mathematically natural true statement unprovable in $\mathbf{PA}$:  \emph{Paris-Harrington Principle} (${\sf PH}$).
Following $\sf PH$, many other mathematically natural statements independent of $\mathbf{PA}$ with combinatorial or number-theoretic contents were formulated: the Kanamori-McAloon principle [Kanamori and McAloon, 1987], the Kirby-Paris sentence [Kirby and Paris, 1982], the
Hercules-Hydra game [Kirby and Paris, 1982], the Worm principle [Beklemishev, 2003][Hamano and Okada, 1997],
the flipping principle [Kirby, 1982], the arboreal statement [Mills, 1980],
P.Pudl\'{a}k's Principle [Pudl\'{a}k, 1979][H\'{a}jek and Paris, 1986], the kiralic and regal principles [Clote and McAloon, 1983] (see [Bovykin, 2006, p. 40]). All these statements are
thought of as much more genuinely and purely mathematical than  G\"{o}del's sentence, and reveal the concrete incompleteness of first order arithmetic.
All these concrete independent arithmetic sentences from mathematics  are provable in fragments of second order arithmetic, and are more complex than G\"{o}del's sentence: G\"{o}del's sentence is equivalent  to $\mathbf{Con}(\mathbf{PA})$ in $\mathbf{PA}$; but all these arithmetic sentences
are not only independent of $\mathbf{PA}$ but also independent of $\mathbf{PA}+ \mathbf{Con}(\mathbf{PA})$ (see [Beklemishev, 2010, p. 36] and [Murawski, 1999, p. 301]).

%In fact, all these principles are equivalent to $\sf PH$ (see [Bovykin, 2006, p.40]).

Harvey Friedman is a leading researcher in the field of  concrete incompleteness. Friedman's work extends the research on concrete incompleteness from first order arithmetic to higher order arithmetic. Friedman's book   ``Boolean Relation Theory and Incompleteness"  is a comprehensive monograph on concrete incompleteness in  mathematics, and provides many examples of concrete mathematical theorems not provable in subsystems of second-order arithmetic stronger than $\mathbf{PA}$, and a number of concrete mathematical statements provable in third-order arithmetic but not provable in second-order arithmetic (see [Friedman, forthcoming]).  Cheng [2015, 2019a]  finds a  concrete mathematical theorem ``Harrington's principle implies the existence of zero sharp", and shows that this theorem  is  expressible in second-order arithmetic, not provable in second-order  or third-order arithmetic, but provable in fourth-order arithmetic.
For more examples of concrete mathematical incompleteness  and  discussions of this  subject, we refer to [Simpson, 1987], [Simpson, 1985], [Pacholski, 1980], [Berline, McAloon and Ressayre, 1981], [Cheng, 2019a], [Bovykin, 2006] and [Friedman, forthcoming].

The impact of G\"{o}del's incompleteness theorem is not confined to the community of mathematicians and logicians, and it has been very popular and widely used outside mathematics and logic.
G\"{o}del's incompleteness theorem has significant philosophical meaning, and raises a number of
philosophical questions concerning the nature of mind and machine, as well as the limit of proof and computation.
G\"{o}del succeeds with his proof of the incompleteness theorem because
he recognizes the central importance of distinguishing theory from metatheory, logic
from meta-logic, signs from their referents (see [Baaz, Papadimitriou, Putnam, Scott and Harper, 2014]). In the literature, there are a lot of discussions about the philosophical meaning of G\"{o}del's incompleteness theorem (see [Gaifman, 2000][Resnik, 1974]
[Auerbach, 1985][Detlefsen, 1979][Detlefsen, 1980][Franks, 2009][Pudl\'{a}k, 1999]). For a popular book about the use and misuse of G\"{o}del's incompleteness theorem in and outside mathematics and logic for a wider audience, we refer to [Franzen, 2005]. In the following, we only give an overview of the philosophical influence of G\"{o}del's incompleteness theorem on the Anti-Mechanism Argument and G\"{o}del's Disjunctive Thesis.

%(4) G\"{o}del's incompleteness theorem has great  influence   on the
%philosophical question of whether the mind can be mechanized;

There are a lot of discussions in the literature about the influence of G\"{o}del's incompleteness theorem on the
philosophical question of whether the mind can be mechanized (see [Penrose, 1989] [Chalmers, 1995] [Lucas, 1996] [Lindstr\"{o}m, 2006] [Feferman, 2009]  [Shapiro, 1998]  [Shapiro, 2003]  [Koellner, 2016] [Koellner, 2018a][Koellner, 2018b] [Cheng, 2020]).
The Anti-Mechanism Argument claims that the mind cannot be mechanized in the sense that the mathematical outputs of the idealized human mind outstrip the mathematical outputs of any Turing machine.\footnote{We will not  consider the
performance of actual human minds, with their limitations and defects; but only consider the idealized human mind and look  at what it can do in principle. See [Koellner, 2018a].}
A popular interpretation of $\sf G1$ is that $\sf G1$ implies that the Anti-Mechanism Argument holds.
G\"{o}del did not argue that his incompleteness theorem implies that the
mind cannot be mechanized.  For G\"{o}del,  the human mind cannot be mechanized and human mind is sufficiently powerful to capture all mathematical truths. G\"{o}del believes that   the distinctiveness of the human mind when compared to a Turing machine is evident in its ability to come up with new axioms and develop new
mathematical theories.  Based on his rationalistic optimism, G\"{o}del believes that we are arithmetically omniscient.
However, G\"{o}del  admits that he cannot give a convincing argument for either the thesis ``the human mind cannot be mechanized" or the thesis ``there are absolutely undecidable statements". G\"{o}del thinks that the most he could claim from his incompleteness theorem is a weaker conclusion, G\"{o}del's Disjunctive Thesis ($\sf GD$),\footnote{The original version of $\sf GD$ was introduce by G\"{o}del in [G\"{o}del, 1995], p. 310: ``So the following disjunctive conclusion is inevitable: either mathematics is incompletable in this sense, that its evident axioms can
never be comprised in a finite rule, that is to say, the human mind
(even within the realm of pure mathematics) infinitely surpasses
the powers of any finite machine, or else there exist absolutely
unsolvable diophantine problems of the type specified (where the
case that both terms of the disjunction are true is not excluded,
so that there are, strictly speaking, three alternatives)".}  which claims that if the human mind can be mechanized, then there are absolutely undecidable statements in the sense that there are
mathematical truths that cannot be proved by the idealized human mind.\footnote{G\"{o}del's Disjunctive Thesis  concerns the limit of mathematical knowledge and the possibility of
the existence of mathematical truths that are inaccessible to the  idealized human mind.} We refer to [Horsten and Welch, 2016] for more discussions of $\sf GD$.

For G\"{o}del, $\sf GD$  is a mathematically established fact of great philosophical interest which follows from his incompleteness theorem, and it is entirely independent from the standpoint taken toward the foundation of mathematics (see [G\"{o}del, 1951]).\footnote{In the literature, there is a consensus that G\"{o}del's argument for $\sf GD$ is definitive, but until now we have
no compelling evidence for or against any of the two disjuncts (see [Horsten and Welch, 2016]).}
For more detailed discussions of G\"{o}del's Disjunctive Thesis and the relationship between G\"{o}del's incompleteness theorem and the Anti-Mechanism Argument, we refer to [Cheng, 2020], [Horsten and Welch, 2016] and Koellner's recent nice papers [Koellner, 2016, 2018a, 2018b].

%From \cite{Papadimitriou}, during the mid-1930s, at least four vastly different definitions of computation were proposed and used to
%demonstrate the impossibility of finding a mechanical procedure for distinguishing between true and false sentences in mathematics: Alan Turing's universal machine; Alonzo Church's
%$\lambda$-calculus; Stephen Kleene's  partial recursive functions; and Emil Post's
%machines.

Now, we give a brief account of  the impact of G\"{o}del's incompleteness theorem on theoretic computer science.
Theoretical computer science is about the power and limitation of
computation.
G\"{o}del's   incompleteness theorem contains several  technical ideas that can
be  recognized as computational. In G\"{o}del's   proof of the incompleteness theorem,
G\"{o}del  uses primitive recursive function and arithmetization which are important tools in theoretic computer science.
The technique of arithmetization which represents syntactic elements, such as logical terms, formulas, and proofs, as
numbers,  has  been used crucially in theoretical
computer science  (see [Papadimitriou, 2014]).
Negative results constitute
an important and distinguishing tradition in theoretical computer science.
One typical example of negative results in theoretic computer science is the undecidability of  the halting problem (the problem of
telling whether a given program will eventually terminate) proved by  Alan Turing.
G\"{o}del's   incompleteness theorem is an ideal archetype, and Turing's halting problem can be seen as a sharpening of G\"{o}del's theorem. For more discussions of the influence of G\"{o}del's   incompleteness theorem on theoretic computer science, we refer to [Baaz, Papadimitriou, Putnam, Scott and Harper, 2014].

\subsection{Fruitfulness}\label{fruitfulness}

%exhibits organizational or explanatory power;
In this section, we discuss another criteria of the depth of G\"{o}del's incompleteness theorem: fruitfulness.
The fruitfulness of a theorem measures   the degree
to which a theorem (or a proof of a theorem) leads to yet further  theorems and
proofs (see [Arana, 2015]).
In this section, we account for the fruitfulness of G\"{o}del's incompleteness theorem from the following three indicators: different proofs of the theorem, generalizations of the theorem (how and in what degree the theorem can be generalized or extended), and the boundary (or the limit) of the theorem (i.e.~ under what conditions the  theorem holds and under what conditions the theorem fails).

The first indicator of the fruitfulness of G\"{o}del's incompleteness theorem is the diversity of its proof methods.
After G\"{o}del, people have found many different proofs of G\"{o}del's incompleteness theorem.
We first give some definitions.
We say a proof of $\sf G1$ is constructive if it explicitly constructs the independent
sentence from the base theory  algorithmically. A non-constructive proof of $\sf G1$ proves the mere existence of the independent sentence, and does not show its existence algorithmically.
We say that a proof of $\sf G1$  has the Rosser property if the proof only assumes that the base theory is consistent instead of assuming that the base theory is  $\omega$-consistent or 1-consistent.

We could classify different proofs of G\"{o}del's incompleteness theorem from the literature based on the following features: (1)  proof-theoretic proof; (2) recursion-theoretic proof;\footnote{For example, Kleene gives a simple proof of $\sf G1$ via recursion theory: for any consistent recursive enumerable theory $T$ that contains $\mathbf{Q}$, there exists some $t \in \omega$ such that $\varphi_t(t)\uparrow$ holds but $T \nvdash ``\varphi_t(t)\uparrow"$ (see  [Salehi and Seraji, 2018, Theorem 2.2]).} (3) model-theoretic proof;\footnote{Arithmetic completeness theorem ([Lindstr\"{o}m, 1997]) is an important tool in the  model-theoretic proof of the incompleteness theorem.} (4) proof via arithmetization; (5) proof via the Diagonalization Lemma;\footnote{Let $T$ be a consistent r.e.~ extension of $\mathbf{Q}$. The Diagonalization Lemma says that for any formula $\phi(x)$ with exactly one free variable, there exists a sentence $\theta$ such that $T\vdash\theta\leftrightarrow\phi(\ulcorner\theta\urcorner)$.} (6) proof based on logical paradox;\footnote{Many  paradoxes have been used to give new proofs of incompleteness theorems: e.g.~ the Liar Paradox, Berry's Paradox, Grelling-Nelson's Paradox and Yablo's Paradox.} (7) constructive proof; (8) proof with the Rosser property; (9) proof via Kolmogorov complexity;\footnote{Kolmogorov complexity is a measure of the quantity of information in finite
objects. Chaitin [1974] gives information-theoretic
formulation of $\sf G1$,  and proves  a weaker version of $\sf G1$ in terms of Kolmogorov complexity.
Kikuchi [1997] proves the formalized version of ${\sf G1}$ via Kolmogorov complexity.} (10) concrete incompleteness (i.e.~ proof via an independent sentence with real mathematical contents).
G\"{o}del's proof of $\sf G1$ has the following features: (1) uses proof-theoretic method with arithmetization; (2) does not directly use the Diagonalization Lemma; (3) the proof formalizes the liar paradox; (4) the proof is constructive;\footnote{I.e.~ given a consistent r.e.~ extension $T$ of $\mathbf{PA}$, one can effectively find a true $\Pi^0_1$ sentence $G_T$ of arithmetic such that $G_T$ is independent of $T$. G\"{o}del calls this the ``incompletability or inexhaustability of mathematics".} (5)  G\"{o}del's   proof does not have the Rosser property; (6) G\"{o}del's sentence is constructed via metamathematical method, and has no real mathematical content.\footnote{I.e.~ G\"{o}del's sentence is a pure logical construction (via the arithmetization of syntax and provability predicate) and has no relevance with classic mathematics (without any combinatorial or number-theoretic content). On the contrary, Paris-Harrington Principle is an independent arithmetic sentence from classic mathematics with combinatorial contents.}

We give some comments about these features of proofs of G\"{o}del's incompleteness theorem. Firstly, these features are not exclusive: a proof of G\"{o}del's theorem may have several above features.
Secondly, each of the above features of G\"{o}del's incompleteness theorem is not a necessary condition to prove  G\"{o}del's theorem.
We have examples of proofs of $\sf G1$ with the above features  and examples of proofs of $\sf G1$ without the above features from the literature.
For example, for the proof of $\sf G1$, we  also have examples of proofs which are non-constructive,\footnote{A non-constructive proof of $\sf G1$ proves the mere existence of the independent sentence and does not show its existence algorithmically.}  proofs having the Rosser property, and proofs  without the  use of arithmetization.\footnote{All  proofs of ${\sf G1}$ we have discussed  use arithmetization. However, A. Grzegorczyk proposes the theory $\mathbf{TC}$ in [Grzegorczyk, 2005] as a possible alternative theory for studying incompleteness
and undecidability, and proves that $\mathbf{TC}$ is incomplete without the use of arithmetization.}
Thirdly, these different proofs of G\"{o}del's incompleteness theorem establish the connection among different fields: proof theory, recursion theory, logical paradox, model theory, Kolmogorov complexity, ect.

The second indicator of the fruitfulness of G\"{o}del's incompleteness theorem is the great variety of its  generalizations.  From the literature, $\sf G1$ and $\sf G2$ can be generalized to both extensions of $\mathbf{PA}$ and weaker theories than  $\mathbf{PA}$ w.r.t.~ interpretation. These generalizations show the applicability and explanatory power of  G\"{o}del's incompleteness theorem. In the following, we give some typical examples to explain this.

We first discuss generalizations of $\sf G1$.
The first example is Rosser's improvement of  $\sf G1$. G\"{o}del's proof of $\sf G1$ assumes that the base theory is $\omega$-consistent.  Rosser proves $\sf G1$ only assuming that the base theory is consistent: Rosser constructs  a Rosser sentence ($\Pi^0_1$ sentence), and shows that if  $T$ is a recursively axiomatized consistent extension of $\mathbf{Q}$, then  the Rosser sentence is independent of $T$. Note that $\omega$-consistency implies consistency. But the converse does not hold, and
the notion of $\omega$-consistency is stronger than consistency since we can find examples of theories that are consistent but not $\omega$-consistent.\footnote{For example, assuming $\mathbf{PA}$ is consistent, then $\mathbf{PA} + \neg \mathbf{Con(PA)}$ is consistent, but not $\omega$-consistent.}

The second example is the generalization of $\sf G1$ to arithmetically definable theories.
From $\sf G1$, if a theory $T$ is a $\Sigma^0_{1}$-definable and consistent extension of $\PA$, then $T$ is not $\Pi^0_{1}$-decisive. Kikuchi-Kurahashi and Salehi-Seraji generalize $\sf G1$ to arithmetically definable theories (see [Kikuchi and Kurahashi, 2017] and [Salehi and Seraji, 2017]), and show that if $T$ is a $\Sigma^0_{n+1}$-definable and $\Sigma^0_{n}$-sound extension of $\mathbf{Q}$, then $T$ is not $\Pi^0_{n+1}$-decisive.

Thirdly, $\sf G1$ can also be generalized via the notion of interpretation. We define that $\sf G1$ holds for a theory $T$ iff for any recursively axiomatizable consistent theory $S$, if $T$ is interpretable in $S$, then  $S$ is incomplete ([Cheng, 2019b]).
In fact, $\sf G1$  also holds for
many weaker theories than  $\mathbf{PA}$ w.r.t.~ interpretation.
Let $\mathbf{R}$ be the theory consisting of schemes $\sf{Ax1}$-$\sf{Ax5}$ with $L(\mathbf{R})=\{\mathbf{0}, \mathbf{S}, +, \times, \leq\}$ where  $\leq$ is a primitive binary relation  symbol, and  $\overline{n}=\mathbf{S}^n \mathbf{0}$ for $n \in \mathbb{N}$:
\begin{description}
  \item[\sf{Ax1}] $\overline{m}+\overline{n}=\overline{m+n}$;
  \item[\sf{Ax2}] $\overline{m}\times\overline{n}=\overline{m\cdot n}$;
  \item[\sf{Ax3}] $\overline{m}\neq\overline{n}$, if $m\neq n$;
  \item[\sf{Ax4}] $\forall x(x\leq \overline{n}\rightarrow x=\overline{0}\vee \cdots \vee x=\overline{n})$;
  \item[\sf{Ax5}] $\forall x(x\leq \overline{n}\vee \overline{n}\leq x)$.
\end{description}

%$\mathbf{Q}$ is finitely axiomatizable but
%$\mathbf{R}$ is not finitely axiomatizable. $\mathbf{Q}$ is essentially incomplete and in fact minimal essentially %incomplete (see \cite{metamathematics}, p.260).

It is well known that $\sf G1$ holds for  $\mathbf{Q}$ and $\mathbf{R}$ (see  [Vaught, 1962]).  For more examples of weaker theories than $\mathbf{PA}$ w.r.t.~ interpretation for which $\sf G1$ holds, we refer to [Cheng, 2019b] for more discussions.

%We say that $\sf G2$ holds for theory $T$ if  is not provable in $T$.

Now, we discuss generalizations of $\sf G2$.  Let $T$ be a recursively axiomatizable consistent extension of $\mathbf{Q}$. Recall that $\mathbf{Con}(T)$ is the canonical arithmetic sentence expressing the
consistency of $T$. In fact, $\sf G2$ can also be generalized in different ways.
Here, we only give two typical examples we think important.
Firstly, $\sf G2$  can be generalized via the notion of interpretation: there is no r.e.~ theory $T$ such that $\mathbf{Q}+\mathbf{Con}(T)$ is interpretable in $T$, i.e.~ $\mathbf{Q}+\mathbf{Con}(T)\ntrianglelefteq T$ (see [Visser, 2011]).
As a corollary, $\sf G2$ holds for any consistent r.e.~  theory  interpreting $\mathbf{Q}$.
Secondly, L\"{o}b's theorem is an important generalization of $\sf G2$.  L\"{o}b shows that for any standard provability predicate $\mathbf{Pr}_T(x)$ and any formula $\phi$, if $T\vdash \mathbf{Pr}_T(\phi)\rightarrow\phi$, then $T\vdash\phi$. As a corollary, we have $T\nvdash \mathbf{Con}(T)$.

%\subsection{Intensionability}\label{on intensionability}

The third indicator of the fruitiness of G\"{o}del's incompleteness theorem is the boundary (or the limit) of the theorem.
The research on the boundary of G\"{o}del's incompleteness theorem reveals the limit of the applicability of G\"{o}del's incompleteness theorem, greatly deepens our understanding of the scope of G\"{o}del's incompleteness theorem, and contributes to new mathematical evidences of the fruitfulness of  G\"{o}del's incompleteness theorem.

We first give a brief account of the boundary (or the limit) of $\sf G1$.
There are many consistent formal theories which are complete.\footnote{For example, the following theories are complete: the theory of dense linear orderings without endpoints ($\mathbf{DLO}$), the theory of ordered divisible groups ($\mathbf{ODG}$), the theory of algebraically closed fields of given characteristic ($\mathbf{ACF_p}$), and the theory of real  closed fields ($\mathbf{RCF}$), etc (see [Epstein, 2011] for details of these theories).}
Whether a theory about arithmetic is complete depends on the language of the theory.  The theory $\PA$ is incomplete in the language  $L(\mathbf{0}, \mathbf{S}, +, \times)$. There are respectively recursively axiomatized complete arithmetic theories in the language of $L(\mathbf{0}, \mathbf{S})$, $L(\mathbf{0}, \mathbf{S}, <)$  and  $L(\mathbf{0}, \mathbf{S}, <, +)$ (see Section 3.1-3.2 in [Enderton, 2001]). Firstly, containing enough information of arithmetic is  essential for the proof of $\sf G1$.\footnote{For example, the Euclidean geometry is not about arithmetic but only about points, circles and lines in general; but the Euclidean geometry is complete as Tarski has proved.}
Secondly, containing the information about the arithmetic of multiplication is essential for the proof of $\sf G1$. If the theory  contains only the information about the arithmetic of addition without multiplication, then it could be complete.\footnote{For example, Presburger arithmetic is the theory of arithmetic of addition, and its language only contains non-logical symbols $\mathbf{0}, \mathbf{S}$ and $+$;
but Presburger arithmetic is complete (see [Murawski, 1999, Theorem 3.2.2]).}  Finally, containing the arithmetic of multiplication is not a sufficient condition for a theory to be incomplete.\footnote{For example, there exists a complete recursively axiomatized theory in the language  $L(\mathbf{0}, \times)$ (see [Murawski, 1999, p.230]).}

Recall that $\sf G1$ holds for some arithmetically definable extensions of $\mathbf{Q}$, but it is not true that any arithmetically definable extension of $\mathbf{Q}$ is incomplete.\footnote{For example, there exists a $\Sigma^0_{n+1}$-definable, $\Sigma^0_{n-1}$-sound  ($n\geq 1$) theory that is a  complete extension of $\mathbf{Q}$ (see  [Salehi and  Seraji, 2017, Theorem 2.6]).}
It was often thought that $\mathbf{R}$ is the weakest theory w.r.t.~  interpretation for which $\sf G1$  holds.
In fact, we can find many theories $S$ weaker than $\mathbf{R}$ such that $\sf G1$  holds for $S$.\footnote{For example, Cheng [2019b] shows that for any recursively inseparable pair $\langle A,B\rangle$, there is a theory $U_{\langle A,B\rangle}$ such that $\sf G1$ holds for $U_{\langle A, B\rangle}$  and $U_{\langle A,B\rangle}\lhd \mathbf{R}$.} We conjecture that there is no minimal r.e.~ theory w.r.t.~ interpretation for which $\sf G1$  holds.

Now, we give a brief account of the boundary (or the limit) of $\sf G2$.
Both mathematically and philosophically, $\sf G2$ is  more problematic than $\sf G1$.
The difference between $\sf G1$ and $\sf G2$ is that, in the case of $\sf G1$, we are mainly interested in the fact that it shows that \emph{some} sentence is independent of the base theory. We make no claim to the effect that that sentence ``really" expresses what we would express by saying ``$\mathbf{PA}$ cannot prove this sentence". But in the case of $\sf G2$, we are also interested in the content of the statement.
In the following, we give a brief overview of the intensionality  of $\sf G2$ (we refer to [Cheng, 2019c] for more details).

For a consistent theory $T$, we say that $\sf G2$ holds for $T$ if the consistency statement of $T$ is not provable in $T$. However, this definition is vague, and whether $\sf G2$ holds for $T$ depends on how we formulate the consistency statement. We refer to this phenomenon as the intensionality of $\sf G2$.
The status of $\sf G2$ is essentially different from $\sf G1$ due to the intensionality of $\sf G2$.
We can say that ${\sf G1}$ is extensional in the sense that we can construct a concrete independent mathematical statement without referring to  arithmetization and provability predicate.
However, ${\sf G2}$ is intensional, and ``whether ${\sf G2}$ holds for $T$" depends on varied factors as we will discuss.

%$\sf G2$ holds under some conditions and fails under some other conditions.
 %However, the meaning of $\mathbf{Con}(T)$  as well as Theorem \ref{G2} and Theorem \ref{general G2}  are vague.

The intensionality of $\sf G2$ has been widely discussed from the literature (e.g.~  [Halbach and Visser, 2014a], [Halbach and Visser, 2014b], [Visser, 2011]).
Visser [2011] locates three
sources of indeterminacy in the formalisation of a consistency statement for a
theory $T$:
(I) the choice of a proof system;
(II) the choice of a way of numbering;
(III) the choice of a specific formula numerating the axiom set of $T$.

In this section, unless stated otherwise, we make the following assumptions:
\begin{itemize}
  \item The theory $T$ is a recursively axiomatized consistent  extension of  $\mathbf{Q}$;
  \item The canonical arithmetic formula to express the consistency of the base theory $T$ is $\mathbf{Con}(T)\triangleq \neg \mathbf{Pr}_T(\mathbf{0}\neq \mathbf{0})$;
  %\item The base proof system is Hilbert-style system with cut elimination;
  \item The canonical numbering we use is G\"{o}del's numbering;
  \item The provability predicate we use is standard;
  \item The formula representing the set of axioms is $\Sigma^0_1$.
\end{itemize}

Based on the current research  on incompleteness from the literature, we argue that ``whether $\sf G2$ holds for $T$" depends on the following factors:
\begin{enumerate}[(1)]
\item  the choice of the base theory $T$;

\item the choice of the method to express consistency;
\item the choice of a provability predicate;

      \item the choice of a numbering;
      \item the choice of a specific formula numerating the axiom set of $T$.
\end{enumerate}

These factors are not independent, and
a choice made at an earlier stage may have effects on the choices
made at a later stage. In the following, when we discuss how  $\sf G2$ depends on one factor,  we always assume that other factors are fixed as in the above default assumptions,   and only the factor we are discussing is varied.
For example,  Visser [2011]  rests on fixed choices for (1)-(2) and (4)-(5) but varies the choice of (3); Grabmayr [2019] rests on fixed choices for (1) and (3)-(5) but varies the choice of (2); Feferman [1960] rests on fixed choices for (1)-(4)   but varies the choice of (5).
In the following, we give a brief discussion of how $\sf G2$ depends on the above five factors.
For more detailed discussions of these factors, we refer to [Cheng, 2019c].

Firstly, ``Whether $\sf G2$ holds for $T$" depends on the choice of the base theory.
A foundational question about $\sf G2$ is: how much of information about arithmetic is required for the proof of $\sf G2$.
If the base theory does not contain enough information of arithmetic,  then $\sf G2$ may fail in the sense that the consistency statement is provable in the base theory.\footnote{Willard [2006]  explores the generality and boundary-case
exceptions of $\sf G2$ under some base theories.
Willard constructs examples of recursively enumerable  arithmetical theories that couldn't prove the totality of successor function but could prove their own canonical consistency (see [Willard, 2001], [Willard, 2006]).}
Pakhomov [2019] defines a theory $H_{<\omega}$, and shows that it proves its own canonical consistency. Thus, $\sf G2$ fails for the theory $H_{<\omega}$.\footnote{Unlike Willard's theories, $H_{<\omega}$ isn't an arithmetical theory but a theory formulated in the
language of set theory with an additional unary function.}

Secondly, ``Whether $\sf G2$  holds for $T$" depends on the choice of the method to express consistency.
From the philosophical point of view, one can ask: what a consistency statement of a theory is?  when can we reasonably say that the arithmetic sentence $\mathbf{Con}(T)$  does really express the consistency of $T$? (see [Visser, 2011, p. 545]).  These questions are difficult to answer, and have been investigated by many logicians, among them Resnik [1974], Detlefsen [1980], Visser [2016, 2011], Feferman [1960], Auerbach [1985] and Franks [2009].
In the literature, we usually use an arithmetic formula in the language of $T$ to express the consistency of $T$.
Artemov [2019] argues  that in Hilbert's consistency program, the original formulation of consistency ``no sequence of formulas  is a derivation of a contradiction" is about finite sequences
of formulas, not about arithmetization, proof codes, and internalized quantifiers.
Artemov concludes that $\sf G2$ does not actually exclude finitary consistency proofs of the original formulation of consistency. Artemov shows that the original formulation of consistency admits a direct proof in informal arithmetic, and this proof is formalizable  in $\mathbf{PA}$ (see [Artemov, 2019]).\footnote{Informal arithmetic is the theory of informal elementary number theory containing recursive identities of addition and multiplication as well as the induction principle. The formal arithmetic $\mathbf{PA}$ is just the conventional formalization of the informal arithmetic (see [Artemov, 2019]).}

In the following, we use a single arithmetic sentence to express the consistency statement.
Even among consistency statements defined via a single arithmetic sentence, we still have different ways to express the consistency of  $T$. For example, another way to express the consistency of  $T$ is
$\mathbf{Con}^0(T) \triangleq \forall x(\mathbf{Fml}(x) \wedge \mathbf{Pr}_T(x) \rightarrow \neg \mathbf{Pr}_T(\dot{\neg} x))$.\footnote{$\mathbf{Fml}(x)$ is the formula which represents the relation that $x$ is a code of a formula.}
Kurahashi [2019] constructs a Rosser provability predicate such that $\sf G2$ holds for the consistency statement formulated via $\mathbf{Con}^0(T)$,\footnote{I.e.~ the consistency statement formulated via $\mathbf{Con}^0(T)$ and the Rosser provability predicate is not provable in $T$.} but $\sf G2$ fails for the consistency statement formulated via $\mathbf{Con}(T)$ (i.e.~ the consistency statement formulated via $\mathbf{Con}(T)$ and the Rosser provability predicate is  provable in $T$).

Thirdly, ``whether $\sf G2$ holds for $T$" depends on the choice of the  provability predicate.
Visser [2016] argues that, being a consistency statement is not an absolute concept but a role w.r.t.~ a choice of the provability predicate (see [Visser, 2016]).
Recall that $\sf G2$ holds for standard provability predicates.\footnote{If a provability predicate $\mathbf{Pr}_{T}(x)$ is standard, then $T\nvdash \neg \mathbf{Pr}_T(\ulcorner \mathbf{0}\neq \mathbf{0}\urcorner)$.}
However, $\sf G2$ may fail for non-standard provability predicates.
Define the consistency statement $\mathbf{Con}^R(T)$ via Rosser provability predicate as $\neg \mathbf{Pr}^R_{T}(\ulcorner \mathbf{0}\neq \mathbf{0}\urcorner)$. Then $\sf G2$ fails for Rosser provability predicate in the sense that  $T\vdash\mathbf{Con}^R(T)$.

Fourthly, ``Whether $\sf G2$ holds for $T$" depends on the choice of numberings.
Any injective function $\gamma$
from a set of $L(\mathbf{PA})$-expressions to $\omega$ qualifies as a numbering.
G\"{o}del's numbering is a special kind of numberings under which the G\"{o}del number of the set of axioms of $\mathbf{PA}$ is recursive.
Grabmayr [2019] shows that $\sf G2$ holds for acceptable numberings.\footnote{For the definition of acceptable numberings, we refer to  [Grabmayr, 2019].} But $\sf G2$  fails for some non-acceptable numberings.

Finally, ``Whether $\sf G2$ holds for $T$" depends on the numeration of $T$. We say that a formula $\alpha(x)$ is a numeration of $T$ if for any $n$, we have $\mathbf{PA} \vdash \alpha(\overline{n})$ iff $n$ is the G\"{o}del number of some $\phi\in T$.
As a generalization, $\sf G2$ holds for any $\Sigma^0_1$ numeration of $T$.\footnote{Given a formula $\alpha(x)$ in $L(T)$,
define the formula $\mathbf{Prf}_{\alpha}(x,y)$ saying ``$y$ is the G\"{o}del number of a proof of the formula with G\"{o}del number $x$ from the set of all sentences satisfying $\alpha(x)$"; define the provability predicate $\mathbf{Pr}_{\alpha}(x)$ of $\alpha(x)$ as $\exists y \mathbf{Prf}_{\alpha}(x,y)$ and the consistency statement $\mathbf{Con}_{\alpha}(T)$ as $\neg \mathbf{Pr}_{\alpha}(\ulcorner \mathbf{0}\neq \mathbf{0}\urcorner)$. This  generalization says that if $\alpha(x)$ is a $\Sigma^0_1$  numeration of $T$, then $T\nvdash \mathbf{Con}_{\alpha}(T)$.}
 However, $\sf G2$ fails for some $\Pi^0_1$ numerations of $T$. For example, Feferman [1960] constructs a $\Pi^0_{1}$ numeration $\tau(u)$ of $T$ such that $\sf G2$ fails under this numeration (i.e.~ $T\vdash \mathbf{Con}_{\tau}(T)$).

\subsection{Unity}

In this section, we discuss the third criteria of the depth of G\"{o}del's incompleteness theorem we propose: unity. The unity of G\"{o}del's incompleteness theorem means that it  ties together apparently disparate fields and draws interconnections between these fields.
In this section, we give a brief account of the unity of G\"{o}del's incompleteness theorem from the following aspects:
\begin{itemize}
  \item finding the order in chaos between mathematics and meta-mathematics;
  \item the close relationship with the theory of undecidability;
  \item the close relationship with logical paradox;
  \item the close relationship with provability logic;
  \item  the close relationship with the formal theory of truth.
\end{itemize}

%establishes the relationship between them

%because its different proofs share features that make manifest connections between different areas of logic and mathematics. Unity

%G\"{o}del's incompleteness theorem is deep because its different proofs share features
%that make manifest connections between different areas of mathematics and logic.

G\"{o}del's proof of
the incompleteness theorem uses  methods from both mathematics and logic. For example, in G\"{o}del's proof, he uses the Chinese Remainder Theorem and the unique factorization in number theory, as well as some meta-mathematical methods in logic such as arithmetization, representability, and self-reference construction.
In Section \ref{section on influence}, we give some examples of natural independent sentences with real mathematical contents.\footnote{For example, the Paris-Harrington Principle, the Kanamori-McAloon principle, the Kirby-Paris sentence, the
Hercules-Hydra game, the Worm principle,
the flipping principle, the arboreal statement,
the P.Pudl\'{a}k's Principle, the kiralic and regal principles.}   Section \ref{section on influence} emphasizes the difference between meta-mathematical independent sentences constructed via pure logic and concrete independent sentences with real mathematical contents.
An interesting and amazing fact is that  all the  mathematically natural independent sentences with  combinatorial or number-theoretic contents we list in Section 3.1 are in fact provably equivalent in $\mathbf{PA}$ to a certain meta-mathematical sentence. Consider the following reflection principle for $\Sigma^0_1$ sentences:
for any $\Sigma^0_1$ sentence $\phi$ in $L(\mathbf{PA})$, if $\phi$ is provable in $\mathbf{PA}$, then $\phi$ is true.
Using the arithmetization of syntax, one can write this principle
as a sentence of $L(\PA)$ and denote it by $\sf Rfn_{\Sigma^0_1}(\mathbf{PA})$ (see [Murawski, 1999, p.301]). McAloon has shown  in $\mathbf{PA}$ that the Paris-Harrington Principle is equivalent to ${\sf Rfn}_{\Sigma^0_1}(\mathbf{PA})$ (see [Murawski, 1999, p. 301]). In fact, similar equivalences can be established for all the   natural independent sentences with combinatorial or number-theoretic contents we list in Section 3.1 (see [Beklemishev, 2010, p. 36], [Beklemishev, 2003, p. 3] and [Murawski, 1999, p. 301]). Isaacson [1987] argues that this result reveals something of the implicit (hidden) higher-order content of the  mathematically natural independent sentences we list in Section 3.1.
This phenomenon shows that the difference between mathematical and meta-mathematical statements is not as huge as we might have expected (see [Dean, 2015] for more discussions of this claim).

In the rest of  this section, we give a brief account of the close relationship among G\"{o}del's incompleteness theorem, the theory of undecidability, logical paradox, provability logic, and the formal theory of truth.

G\"{o}del's work has  many crucial connections to the theory of computation and undecidability.
G\"{o}del's proof  contains the germs of such influential computational ideas as arithmetization and primitive recursion.
The method of arithmetization plays a major role in the growth of recursion theory. We say a theory $T$ is
\emph{essentially undecidable} if any recursively axiomatizable consistent extension of $T$ in the same language is undecidable; and  $T$ is \emph{essentially incomplete}  if  any recursively axiomatizable consistent extension of $T$ is incomplete.
Since a theory $T$ is essentially undecidable  if and only if $T$ is essentially incomplete,\footnote{This follows from the following two facts: (1) every consistent recursively axiomatizable
complete theory is decidable; (2) every incomplete decidable theory has a consistent,
decidable complete extension in the same language (see  [Murawski, 1999, p. 214-215]).} the theory of completeness/incompleteness  is closely related to the theory of decidability/undecidability.
Recall that we have defined the notion ``$\sf G1$ holds for a theory $T$".
From [Cheng, 2019b], $\sf G1$ holds for  $T$
if and only if $T$ is essentially undecidable.
It is well known that we can prove $\sf G1$ and $\sf G2$ in terms of the undecidability of the halting problem. All of these show the close relationship between the theory of incompleteness and the theory of undecidability.

The current research practice on G\"{o}del's incompleteness theorem reveals that $\sf G1$ is closely related to logical paradox.
G\"{o}del comments in his famous paper  that ``any epistemological antinomy could be used for a similar proof of the existence of undecidable propositions" (see [Feferman, 1995a]).
In G\"{o}del's proof of ${\sf G1}$, we can view G\"{o}del's sentence as the formalization of the Liar Paradox.
G\"{o}del's sentence concerns the notion of provability, but the liar sentence in the Liar Paradox concerns the notion of truth in the standard model of arithmetic.
Except for the Liar Paradox, many other paradoxes have been properly formalized to give new proofs of the incompleteness theorem: for example, Berry's Paradox in [Boolos, 1989][Chaitin, 1974][Kikuchi, 1994] [Kikuchi, Kurahashi and Sakai, 2012] [Kikuchi and Tanaka, 1994][Vopenka, 1966], Grelling-Nelson's Paradox in [Cie\'{s}li\'{n}ski, 2002], the Unexpected Examination Paradox in [Fitch, 1964] [Kritchman and Raz, 2010], and Yablo's Paradox in [Cie\'{s}li\'{n}ski and Urbaniak, 2013] [Kurahashi, 2014][Kurahashi, 2014][Priest, 1997].
%Dean [2019] gives a detailed discussion on how to use the arithmetized completeness theorem to give formal incompleteness results from set theoretic and semantic paradoxes.

%Incompleteness theorem and formal theory of truth.

One important consequence of G\"{o}del's incompleteness theorem  is Tarski's undefinability theorem of truth as an application of the Diagonalization Lemma. Define $\textbf{Prov}=\{\phi\in L(\PA): \mathbf{PA}\vdash\phi\}$ and $\textbf{Truth}=\{\phi\in L(\PA):\mathfrak{N}\models\phi\}$. From Tarski's theorem, \textbf{Truth} (the set of true sentences of arithmetic) is not definable in the standard model of arithmetic; as a corollary, \textbf{Truth} is not arithmetic and not representable in $\mathbf{PA}$. But \textbf{Prov} (the set of sentences provable in $\mathbf{PA}$) is definable in the standard model of arithmetic and recursive enumerable, even if it is  not recursive  (for details of  properties of \textbf{Truth} and \textbf{Prov}, we refer to [Murawski, 1999] [Tarski,  Mostowski and Robinson, 1953]).
Current research practice reveals  the relationship between G\"{o}del's incompleteness theorem  and Tarski's undefinability theorem of truth. For example, Visser [2019] gives a self-reference-free proof of G\"{o}del's second incompleteness theorem  from  Tarski's undefinability theorem of truth.

%Incompleteness theorem and provability logic.

Provability logic is an important tool for the study of incompleteness and meta-mathematics of arithmetic.
The origins of provability logic (e.g.~ Henkin's problem, the isolation of derivability conditions, L\"{o}b's theorem) are all closely tied to the incompleteness theorem historically. In this sense, we can say that G\"{o}del's incompleteness theorem plays a unifying role between first order arithmetic and modal logic.
The notion of arithmetical interpretation provides us with an important tool to establish the relationship between provability logic and meta-mathematics of arithmetic.\footnote{Let $T$ be a consistent r.e.~ extension of $\mathbf{Q}$. A mapping from the set of all modal propositional variables to the set of $L(T)$-sentences is called an \emph{arithmetical interpretation}. Every arithmetical interpretation $f$ is uniquely extended to the
mapping $f^{\ast}$ from the set of all modal formulas to the set of $L(T)$-sentences so
that $f^{\ast}$ satisfies the following conditions:
(1) $f^{\ast}(p)=f(p)$ for each propositional variable $p$; (2) $f^{\ast}$ commutes with every propositional connective; (3) $f^{\ast}(\Box A)$ is $\mathbf{Pr}_T(\ulcorner f^{\ast}(A)\urcorner)$ for every modal formula $A$.
We equate an arithmetical interpretation $f$ with its unique extension $f^{\ast}$ defined on the set of all modal formulas.}
Surprisingly, Solovay's Arithmetical Completeness Theorems for $\mathbf{GL}$ and $\mathbf{GLS}$\footnote{The definition of $\mathbf{GL}$ and $\mathbf{GLS}$ is standard, and we refer to [Boolos, 1993].} characterize the difference between $\textbf{Prov}$ and $\textbf{Truth}$ via provability logic.\footnote{Solovay's Arithmetical Completeness Theorem for $\mathbf{GL}$ says that if $T$ is a $\Sigma^0_1$-sound r.e.~ extension of $\mathbf{Q}$, then for any modal formula $\phi$ in $L(\mathbf{GL})$, $\mathbf{GL}\vdash\phi$ iff $T\vdash \phi^{f}$ for every arithmetic interpretation $f$.
Solovay's Arithmetical Completeness Theorem for $\mathbf{GLS}$ says that for any modal formula $\phi$, $\mathbf{GLS}\vdash\phi$ iff $\mathfrak{N}\models\phi^{f}$ for every arithmetic interpretation $f$.}

Provability logic is the logic of properties of provability predicates. Note that the proof of G\"{o}del's incompleteness theorem depends on the property of provability predicates.  Provability logic provides us with a new perspective and an important tool  to understand incompleteness. Provability logics based on different provability predicates reveal the intensionability of provability predicates which is one source of the intensionability of $\sf G2$, and provide us with a new route to examine the intensionability of provability predicates. Under different numerations of the base theory, the provability predicate may have different properties, and hence may correspond to different  provability logics (i.e.~ different modal principles under arithmetic interpretations). For more discussions about the relationship between incompleteness and provability predicate, we refer to [Kurahashi, 2019a] [Kurahashi, 2019b].

\section{Some explanations}

In this section, we give some explanations for our account of the depth of G\"{o}del's incompleteness theorem.

In the literature, some criteria for mathematical depth have been proposed.  From [Ernst, Heis, Maddy, McNulty,
and Weatherall, 2015b], the following five candidate criteria for mathematical depth  have gained the widest
support:
\begin{enumerate}[(1)]
  \item ties together apparently disparate fields;
  \item involves impurity (definitions that reach into higher types, proofs that appeal
to concepts other than those in the statement proved);
  \item finds order in chaos;
  \item exhibits organizational or explanatory power;
  \item transforms a field or opens a new one.
\end{enumerate}

No examples of mathematical depth that failed on all above five criteria were proposed. From our account of the depth of G\"{o}del's incompleteness theorem, G\"{o}del's incompleteness theorem  satisfies all of the above criteria. Compare our three criteria with the above five criteria, our three criteria are more  general than the above five criteria.
Each one of our three criteria is   not a sufficient  condition for a theorem to be deep.
We can not say that if a theorem is influential (or fruitful, or exhibiting unity), then this theorem is deep. For example, a theorem may be fruitful but not widely considered as deep.
It is controversial  whether each one of our three criteria is a necessary condition for a theorem to be deep (i.e.~ is it true that if a theorem is not influential (or not fruitful, or not exhibiting unity), then this theorem is not deep).
We do not know whether there might be an example of depth without fruitfulness.

%A natural question is whether our three criteria of the depth of G\"{o}del's incompleteness theorem can also be accepted as general criteria of mathematical depth.
 %For us, our three criteria (influence, fruitfulness, unity) are jointly sufficient and necessary conditions for a theorem to be deep.

%But the problem is:  we have no well-formulated measure of the degree of the influence, fruitfulness and unity of mathematical theorems and it is hard to compare the depth of different theorems. For example, the influence of G\"{o}del's incompleteness theorem covers logic, mathematics, philosophy and theoretic computer science, and the influence of Szemer\'{e}di's Theorem covers different fields of pure mathematics. How can we say that one is more influential than another one, or one is deeper than another one?

The depth of a theorem is not an essential intrinsic  property of the theorem but a property of mathematical practice of this  theorem.
Depth is historically located or contextual. What people
once thought deep could turn out not to be. For example, the theorem on the irrationality of the square root of two was regarded as deep by the ancients, but perhaps today it appears too simple to be deep (see [Ernst, Heis, Maddy, McNulty,
and Weatherall, 2015a]).
A theorem may not seem deep immediately after its first publication, but may be widely considered as deep during the mathematical practice. Deep theorems are generally the work of several generations of mathematicians. For example, for the influence, fruitfulness  and unity of  G\"{o}del's incompleteness theorem, G\"{o}del may not even realize that his theorem is so influential on foundations of mathematics, philosophy, mathematics and theoretic computer science; has so many different proofs and generalizations; establishes so many connections among varied fields;  and whether $\sf G2$ holds depends on so many factors.
Nowadays, G\"{o}del's incompleteness theorem and its proof are standard materials of  logic textbooks for advanced undergraduates.
The current research practice on incompleteness (such as concrete incompleteness) is rather complex or even more  technical than G\"{o}del's original proof. Thus, if  the depth of G\"{o}del's incompleteness theorem is only linked to the original proof by G\"{o}del, and is not related to the research practice of this theorem after G\"{o}del, then we may no longer view G\"{o}del's incompleteness theorem as deep since the later research practice on incompleteness has greatly deepen our understanding of G\"{o}del's incompleteness  theorem.
Finally, influence, fruitfulness and unity of a theorem all depend on the level of research practice of this theorem. But there is no limit of research practice, and as research practice goes on, people may find more and more evidences of  influence, fruitfulness and unity of this theorem.

%A theorem may not be influential or fruitful or exhibiting unity in our sense immediately  after its publication, but it may become more and more influential, have more distinct proofs, lead to more results and establish more unity during the ongoing research practice.

%Assuming we view an account of mathematical depth as subjective if our justifications depend on the status of the current research practice, then our accounts of Influence-Fruitfulness-Unity of G\"{o}del's incompleteness theorem are subjective. The key point is that

%The subjectivity and objectivity of our four accounts of the depth of G\"{o}del's incompleteness theorem.

A natural question is:  is depth an objective property (independent of our interests and abilities) or a subjective property of the theorem (something essentially tied to our interests, abilities, and so on)? It is not our goal to decide on the larger question of whether
depth is objective and what objectivity  would consist in.
However, it is an interesting question  whether our account of the depth of G\"{o}del's incompleteness theorem is objective (and if so, in what sense), and whether it can deliver a notion of depth that is not essentially dependent on our contingent interests and abilities.
For us, this depends on how we view the objectivity of  our account of the depth of G\"{o}del's incompleteness theorem.
Assuming we view our account as objective if our justifications are based on  mathematical evidences of the theorem, then  our  account of the depth of G\"{o}del's incompleteness theorem is objective since our account  is  based on mathematical evidences from the research practice of   G\"{o}del's theorem,
not based on individual  preferences, interests and abilities (even if we have limited knowledge about the current research of G\"{o}del's incompleteness theorem  and the mathematical evidences we give here  are limited).
The evaluation of the depth of a mathematical theorem may be person-dependent: for the same theorem, some may think it as interesting and deep, but others may  not.
For mathematicians interested in foundations of mathematics, they may view G\"{o}del's incompleteness theorem as influential; but for mathematicians without any interest on foundations of mathematics, they may not view  G\"{o}del's incompleteness theorem as  influential.
But we may have an objective account of mathematical depth of a theorem from some academic community.
The depth of a mathematical theorem is more than a fact about the theorem but an evaluation of this theorem from the specific academic community.  For example, the depth of G\"{o}del's incompleteness theorem
should  not be judged  by the academic community from topology; instead, it should be  judged  by the academic community from mathematical logic.
We can describe a possible practical procedure of judging whether a given mathematical theorem is deep. Given a mathematical theorem $A$ in some  field $X$, whether theorem $A$ is deep can be judged by a group of  academic committees  which consist of  top scholars around the world in the field $X$. Following the commonly accepted criteria of mathematical depth by this group (such as our Influence-Fruitfulness-Unity criteria), this group of academic committees can rank the depth of  theorem $A$ according to the commonly accepted criteria.

%Assuming we view our account of the depth of G\"{o}del's incompleteness theorem as objective if our account does not depend on , then our account of the depth of G\"{o}del's incompleteness theorem is  objective.

%Thus if depth
%is interest-dependent, a theorem is deep only so long as those consequences are of
%interest to us, and this threatens to render depth subjective.

%However, the accounts for the depth of G\"{o}del's incompleteness theorem in this paper are not personal but
%has objective, or at least considerable inter-subjective, validity.

Another natural question is: are there salient differences between the depth of G\"{o}del's incompleteness theorem and the depth of a pure mathematical theorem such as Szemer\'{e}di's Theorem?
Based on Arana's work in [2015], we can argue that Szemer\'{e}di's Theorem also satisfies our Influence-Fruitfulness-Unity criteria.\footnote{This is an observation based on Arana's work in [2015]. Due to the space limit, we can not provide the details of the argument for this claim.}
Thus, according to our Influence-Fruitfulness-Unity criteria, both Szemer\'{e}di's Theorem and G\"{o}del's incompleteness theorem are deep.
However, even if both Szemer\'{e}di's Theorem and G\"{o}del's incompleteness theorem satisfy our Influence-Fruitfulness-Unity criteria, but they have different justifications for the three criteria.
For example, even if both Szemer\'{e}di's Theorem and G\"{o}del's incompleteness theorem  are influential, their influence cover different fields. The impact of Szemer\'{e}di's Theorem is more on pure mathematics and especially number theory. But the impact of G\"{o}del's incompleteness theorem is more on logic. We do not know any general criterion of mathematical theorems which Szemer\'{e}di's Theorem satisfies but G\"{o}del's incompleteness theorem does not satisfy.
%The difference among their justifications reflects the difference between the depth of G\"{o}del's incompleteness theorem and the depth of Szemer\'{e}di's Theorem.

%Since Szemer\'{e}di's Theorem and G\"{o}del's incompleteness theorem are from two different fields in mathematics and they have distinct justifications for their depth, it does not make sense to say that one is deeper than another.

In  summary, in this paper, we put forward G\"{o}del's incompleteness theorem as a case for studying mathematical depth.
We propose three criteria (influence, fruitfulness, and unity) to account for the depth of G\"{o}del's incompleteness theorem based on the current research practice, and justify that G\"{o}del's theorem satisfies our Influence-Fruitfulness-Unity criteria.
Many points discussed in this paper are worth further exploration. For example, the uniform criteria of mathematical depth, the objectivity of mathematical depth, the difference between the depth of theorems and the depth of proofs, the method to compare the depth of different mathematical theorems, and the difference between the depth of pure logical theorems and the depth of pure mathematical theorems are all worthy of further study. Deeper research about these topics might bring more insights of mathematical depth to light.
In this paper, we focus on the methodological study of what the depth of G\"{o}del's incompleteness theorem consists in?
We hope our account of the depth of  G\"{o}del's incompleteness theorem sheds a little light on mathematical depth as a notion with many faces.

%A fuller account will have to come to grips with the plurality of analysis raised here, to see if perhaps some further deeper criterions can be found for the account of mathematical account.

%opens up new lines of inquiry, that facilitates proofs, that provides explanations

%This paper, we focus on the discussion of informal depth.

%Logic has demonstrated
%that there are hierarchies of provability, complexity, and unsolvability.

%There are two levels of depth of a theorem: formal depth and informal depth. Formal depth measures the depth of the theorem via the intrinsic property of the theorem: for example,

\section*{Funding}
This work was supported by The National Social Science Fund of China for general project ``Research on the limit of incompleteness and the intensional problem of incompleteness" [grant number 18BZX131].

\section*{Acknowledgements}
I would like to thank  the referees for providing detailed and  helpful comments for improvements. I would like to thank the fund support by The National Social Science Fund of China for general project.

\section*{Reference}

%\bibitem{Arana}
Arana, Andrew [2015]: `On the Depth of Szemer\'{e}di's Theorem', \emph{Philosophia Mathematica}, Volume 23, Issue 2,  Pages 163-176.

Artemov, Sergei [2019]: `The Provability of Consistency', reprint, see
arXiv:1902.07404v5, 2019.

Auerbach, D.David  [1985]:
`Intensionality and the G\"{o}del Theorems',
\emph{Philosophical Studies: An International Journal for Philosophy in the Analytic
Tradition}, Vol. 48, No. 3, pp. 337-351.

Avigad, Jeremy [2005]:  `Incompleteness via the halting problem', 2005.

Baaz, Matthias,  Christos H. Papadimitriou, Hilary W. Putnam, Dana S. Scott, and Charles L. Harper Jr (Edited)[2014]:  `Kurt G\"{o}del and the Foundations of Mathematics: Horizons of Truth', Cambridge University Press.

Beklemishev, D.  Lev [2003]: `The Worm principle', Logic Group Preprint Series,
219, Utrecht Univ. March.

Beklemishev, D.  Lev [2010]: `G\"{o}del incompleteness theorems and
the limits of their applicability I',
\emph{Russian Math Surveys}.

Berline, C.,  McAloon, K., Ressayre, J.P.  (editors) [1981]: \emph{`Model Theory and Arithmetic'},
Lecture Notes in Mathematics, vol. 890, Springer, Berlin.

%\bibitem{Boolos 90}
%George Boolos et al. An open peer commentary on The Emperor's new
%mind. Behavioral and Brain Sciences, 13(4):655-91, December 1990.

Boolos, George [1989]: `A new proof of the G\"{o}del incompleteness theorem', \emph{Notices Amer. Math.
Soc}.36, 388-390.

Boolos, George [1993]:
\emph{`The Logic of Provability'}, Cambridge University Press.

Bovykin, Andrey [2006]: `Brief introduction to unprovability',
 Logic Colloquium, Lecture Notes in Logic 32.

%Bernd Buldt. The Scope of G\"{o}del's First Incompleteness Theorem.
%\emph{Logica Universalis}. 8 (3-4), 499-552, 2014.

Chaitin, J. Gregory [1974]: `Information-theoretic limitations of formal systems',
\emph{Journal of the Association for Computing Machinery}, 21:403-424.

Chalmers, David J [1995]: `Minds, machines, and mathematics: A review of
Shadows of the mind by Roger Penrose', \emph{Journal Psyche}, 2.

Chao, Conden  and Seraji  Payam [2018]: `G\"{o}del's second incompleteness theorem for $\Sigma_{n}$-definable theories',
\emph{Logic Journal of the IGPL}, Volume 26, Issue 2, Pages 255-257.

%\bibitem{Chang, Chen Chung 90}
%Chang, Chen Chung; Keisler, H. Jerome. Model Theory. Studies in Logic and the Foundations of Mathematics (3rd ed.). Elsevier, 1990.

Cheng, Yong [2019a]: \emph{`Incompleteness for Higher-Order Arithmetic: An Example Based on Harrington's Principle'}, Springer series: Springerbrief in Mathematics, Springer.

Cheng, Yong [2019b]: `Finding the limit of incompleteness I', accepted and to appear in \emph{The Bulletin of Symbolic Logic}, DOI: 10.1017/bsl.2020.09.

Cheng, Yong [2019c]: `Current research on G\"{o}del's incompleteness theorem', accepted and to appear in The Bulletin of Symbolic Logic, DOI: 10.1017/bsl.2020.44.

Cheng, Yong and Schindler Ralf [2015]: `Harrington's Principle in higher order arithmetic', {\em The Journal of Symbolic Logic,} Volume 80,Issue 02, pp 477-489.

 Cie\'{s}li\'{n}ski, Cezary [2002]: `Heterologicality and incompleteness', \emph{Mathematical Logic Quarterly}, 48(1),
105-110.

Cie\'{s}li\'{n}ski, Cezary and  Urbaniak, R [2013]: `Godelizing the Yablo sequence', \emph{Journal of Philosophical Logic}, 42(5), 679-695.

%\bibitem{Davis 93} Martin Davis et al. Continuing commentary on The Emperor's new
%mind. Behavioral and Brain Sciences, 16(3):611-16, Semptermber 1993.

Clote, P. and McAloon, K [1983]: `Two further combinatorial theorems equivalent to the 1-consistency of Peano arithmetic', \emph{J. Symb. Log}. vol.48, no.4, pp. 1090-1104.

Dean, Walter [2019]: `Incompleteness via paradox and completeness', to appear in \emph{Review of symbolic logic}.

Detlefsen, Michael [1979]: `On Interpreting G\"{o}del's Second Theorem',
\emph{Journal of Philosophical Logic} 8, 297-313.

Detlefsen, Michael [1980]: `On a theorem of Feferman',
\emph{Philosophical Studies: An International Journal for Philosophy in the Analytic
Tradition}, Vol. 38, No. 2, pp. 129-140.

Enderton B. Herbert [2001]:
\emph{`A mathematical introduction to logic'} (2nd ed.), Boston, MA: Academic Press.

Epstein L. Richard [2011]  (with contributions by Les{\l}aw W.Szczerba):
\emph{`Classical mathematical logic: The semantic foundations of logic'},  Princeton University Press, 2011.

Ernst, Michael, Heis Jeremy , Maddy Penelope, McNulty  Bennett Michael, Weatherall Owen James [2015a]: `Foreword to Special Issue on Mathematical Depth', \emph{Philosophia Mathematica}, Volume 23, Issue 2, Pages 155-162.

Ernst, Michael, Heis Jeremy , Maddy Penelope, McNulty  Bennett Michael, Weatherall Owen James [2015b]: `Afterword to Special Issue on Mathematical Depth',
\emph{Philosophia Mathematica}, Volume 23, Issue 2,  Pages 242-254.

Feferman, Solomon [1960]: `Arithmetization of mathematics in a general setting',
\emph{Fundamenta Mathematicae}, vol. 49, pp. 35-92.

Feferman, Solomon [1988]: `Hilbert's Program Relativized: Proof-Theoretical and Foundational Reductions',
\emph{The Journal of Symbolic Logic}, Vol. 53, No. 2, pp. 364-384.

Feferman, Solomon [2006]:
`The Impact of the
Incompleteness Theorems
on Mathematics',
\emph{Notices of the AMS}, Volume 53, Number 4, p.434-439.

Feferman, Solomon [2009]: `G\"{o}del, Nagel, minds, and machines',
\emph{The Journal of
Philosophy}, CVI(4):201-219.

Feferman, Solomon,  Dawson, Jr. W. John, Goldfarb Warren, Parsons Charles, and Solovay M. Robert (Edited) [1995a]: 'Kurt G\"{o}del's Collected Works (vol. 1, pp. 145-195)',  Oxford University Press, New York and Oxford.

Feferman, Solomon,  Dawson, Jr. W. John, Goldfarb Warren, Parsons Charles, and Solovay M. Robert (Edited) [1995b]: `Collected Works, Volume III: Unpublished Essays and Lectures', Oxford University Press, New York and Oxford.

Fitch, F.B [1964]: `A goedelized formulation of the prediction paradox', \emph{American Philosophical
Quarterly}, 1, 161-164.

Franks, Curtis [2009]:
\emph{`The Autonomy of Mathematical Knowledge: Hilbert's Program
Revisited'}, Cambridge University Press.

Franzen, Torkel [2005]:
\emph{`G\"{o}del's Theorem: an incomplete guide to its use and abuse'}, A.K.Peters.

Friedman, M. Harvey [formthcoming]:
\emph{`Boolean Relation Theory and Incompleteness'}, Lecture Notes in Logic, Association for Symbolic Logic.

Gaifman, Haim [2000]: `What G\"{o}del's incompletness result does and does not
show',
\emph{The Journal of Philosophy}, XCVII(8):462-470.

G\"{o}del, Kurt [1931]: `\"{U}ber formal unentscheidbare S\"{a}tze der Principia Mathematica und
verwandter Systeme I',
\emph{Monatsh. Math. Phys}. 38:1, 173-198.

G\"{o}del, Kurt [1995]: `Some basic theorems on the foundations of mathematics and their implications', in `Collected Works, Volume III: Unpublished Essays and Lectures', pp. 304-323, Oxford
University Press.

Grabmayr, Balthasar [2019]: `On the Invariance of G\"{o}del's Second Theorem
with regard to Numberings', Reprint.

Gray, Jeremy [2015]:  `Depth-A Gaussian Tradition in Mathematics',
\emph{Philosophia Mathematica}, Volume 23, Issue 2, Pages 177-195.

Grzegorczyk, Andrzej [2005]: `Undecidability without arithmetization', \emph{Studia Logica}, 79(2):163-230.

H\'{a}jek, Petr and Paris J [1986]: `Combinatorial principles concerning approximations of functions',
\emph{Archive Math. Logic}, vol.26, no.1-2, pp.13-28.

H\'{a}jek, Petr and Pudl\'{a}k Pavel [1993]:
\emph{`Metamathematics of First-Order
Arithmetic'}, Springer-Verlag, Berlin-Heidelberg-New York.

Halbach, Volker and Visser Albert [2014a]:   `Self-reference in arithmetic I',
\emph{Review of
Symbolic Logic} 7(4), 671-691.

Halbach, Volker and Visser Albert [2014b]:  `Self-Reference in Arithmetic II', \emph{Review of
Symbolic Logic} 7(4), 692-712.

Hamano, M  and M Okada [1997]:  `A relationship among Gentzen's proof-reduction,
Kirby-Paris' hydra game, and Buchholz's hydra game',
\emph{Math. Logic Quart}, 43:1, 103-120.

Horsten, Leon  and Philip Welch [2016]:
\emph{`G\"{o}del's Disjunction: The scope and limits of mathematical knowledge'},  Oxford University Press.

Isaacson, D. [2011]: `Necessary and sufficient conditions for undecidability of the G\"{o}del sentence and its truth', In: D. DeVidi, etal. (Eds.), Logic, Mathematics, Philosophy: Vintage Enthusiasms, Springer, ISBN 9789400702134, pp. 135-152.

Kanamori, A  and McAloon K.[1987]:  `On G\"{o}del's incompleteness and finite
combinatorics',
\emph{Ann. Pure Appl. Logic}, 33:1, 23-41.

Kikuchi, Makoto [1994]: `A note on Boolos' proof of the incompleteness theorem',  \emph{Math. Log. Quart}. 40, 528-532.

Kikuchi, Makoto [1997]: `Kolmogorov complexity and the second incompleteness
theorem', \emph{Archive for Mathematical Logic}, 36(6):437-443.

Kikuchi, Makoto and Kurahashi Taishi  [2017]: `Generalizations of G\"{o}del's incompleteness theorems for $\Sigma_n$-definable theories of arithmetic',
\emph{Rew. Symb. Logic}, 10(4): 603-616.

Kikuchi, Makoto, Kurahashi Taishi and Sakai, H [2012]:  `On proofs of the incompleteness theorems based on Berry's paradox by Vop\v{e}nka, Chaitin, and Boolos', \emph{Mathematical Logic Quarterly}, 58(4-5), 307-316.

Kikuchi, Makoto and Tanaka Kazuyuki [1994]: `On formalization of model-theoretic proofs of G\"{o}del's theorems', \emph{Notre Dame Journal of Formal Logic},
35(3):403-412.

Kirby, L [1982]: `Flipping properties in arithmetic', \emph{J. Symb. Log}. vol.47,
no.2, pp. 416-422.

Kirby, L. A. S.  and Paris J. [1982]: `Accessible independence results for Peano arithmetic',
\emph{Bull. London Math. Soc}, 14:4, 285-293.

Koellner, Peter [2016]: `G\"{o}del's Disjunction. Charpter in G\"{o}del's Disjunction: The scope and limits of mathematical knowledge', edited by Leon Horsten and Philip Welch, Oxford University Press.

Koellner, Peter [2018]: `On the Question of Whether the Mind Can Be
Mechanized Part I: From G\"{o}del to Penrose',
\emph{Journal of Philosophy},
Volume 115, Issue 7, Pages 337-360.

Koellner, Peter [2018]:  `On the Question of Whether the Mind can be
Mechanized Part II: Penrose's New Argument', \emph{Journal of Philosophy},
Volume 115, Issue 9, Pages 453-484.

Kotlarski, Henryk [2004]:
`The incompleteness theorems after 70 years',
\emph{Annals of Pure and Applied Logic} 126, 125-138, 2004.

Kritchman, S.  and  Raz, R. [2010]: `The surprise examination paradox and the second incompleteness
theorem',
\emph{Notices of the American Mathematical Society}, 57(11), 1454-1458.

Kurahashi, Taishi [2014]:   `Rosser-type undecidable sentences based on Yablo's paradox',
\emph{Journal of
Philosophical Logic}, 43(5), 999-1017.

Kurahashi, Taishi [2019a]:  `A note on derivability conditions', Reprint.

Kurahashi, Taishi [2019b]: `Rosser provability and the second incompleteness
theorem', Reprint.

Lange, Marc [2015]:  `Depth and Explanation in Mathematics',
\emph{Philosophia Mathematica}, Volume 23, Issue 2,  Pages 196-214.

Lindstr\"{o}m, Per [1997]: \emph{`Aspects of Incompleteness'}, Lecture Notes in Logic v. 10.

Lindstr\"{o}m, Per [2006]:  `Remarks on Penrose's new argument',
\emph{Journal of
Philosophical Logic}, 35:231-237.

Lucas, R. J. [1996]: `Minds, machines, and G\"{o}del: A retrospect', in 'Machines and thought: The legacy
 of Alan Turing, Volume 1' (P. J. R. Millican and A. Clark, editors), Oxford University Press,
 Oxford.

Mills, G. [1980]: `A tree analysis of unprovable combinatorial statements', \emph{Model Theory of Algebra
and Arithmetic}, Lecture Notes in Mathematics, vol.834, Springer, Berlin pp.248-311.

Murawski, Roman [1999]:
\emph{`Recursive Functions and Metamathematics: Problems of Completeness and Decidability, G\"{o}del's Theorems'}, Springer Netherlands.

%\bibitem{Theory of recursive functions} {H. Rogers,} Theory of recursive functions and effective computability, McGraw-Hill, New York, 1967.

Pacholski, L. and Wierzejewski,J. [1980]: \emph{`Model Theory of Algebra and Arithmetic'},
Lecture Notes in Mathematics, vol.834, Springer, Berlin.

Pakhomov, Fedor [2019]: `A weak set theory that proves its own consistency', Reprint.

Papadimitriou, H. Christos [2014]: `Computation and Intractability:
Echoes of Kurt Godel', in 'Kurt G\"{o}del and the Foundations of Mathematics: Horizons of Truth'.

Paris, J.  and Harrington, L. [1977]: `A mathematical incompleteness in Peano arithmetic',
Handbook of mathematical logic (J. Barwise, ed.), Stud. Logic Found. Math.,
vol. 90, North-Holland, Amsterdam-New York-Oxford, pp. 1133-1142.

Pudlak, Pavel [1979]: `Another combinatorial principle independent of Peano's axioms', Unpublished.

Pudlak, Pavel [1999]: `A note on applicability of the incompleteness
theorem to human mind',
\emph{Annals of Pure and Applied Logic} 96, 335-342.

Penrose, Roger [1989]:
\emph{`The Emperor's New Mind: Concerning Computeres,
Minds, and the Laws of Physics'}, Oxford University Press.

%\bibitem{Penrose 90}
%Roger Penrose. Pr\'{e}cis of The Emperors new mind-concerning com-
%puters, minds and the laws of physics. Behavioral and Brain Sceinces,
%13(4):643-54, December 1990.

%\bibitem{nonalgorithmic mind} Roger Penrose. The nonalgorithmic mind. Behavioral and Brain
%Sceinces, 13(4):692-706, December 1990.

%\bibitem{Penrose 93} Roger Penrose. An Emperor still without mind. Behavioral and Brain
%Sceinces, 16(3):616-622, September 1993.

Priest, G [1997]: `Yablo's paradox', \emph{Analysis}, 57(4), 236-242.

%\bibitem{Reinhardt 86} {W.N. Reinhardt,} Epistemic theories and the interpretation of \Godel's incompleteness theorems. Journal of Philosophical Logic, 15 (4):427¨C474, 1986.

Resnik, D. M. [1974]:  `On the Philosophical Significance of Consistency Proofs',
\emph{Journal of Philosophical Logic} 3, 133-147.

Rosser, Barkley John [1936]: `Extensions of some theorems of G\"{o}del and
Church', \emph{The Journal of Symbolic Logic}, 1(3):87-91.

Salehi, Saeed  and Seraji Payam [2017]: `G\"{o}del-Rosser's Incompleteness Theorem, generalized and optimized for definable theories',
\emph{Journal of Logic and Computation}, Volume 27, Issue 5, Pages 1391-1397.

Salehi, Saeed  and Seraji Payam [2018]: `On constructivity and the Rosser property: a closer look at some G\"{o}delean proofs', \emph{Annals of Pure and Applied Logic} 169 (2018) 971-980.

Shapiro, Stewart [1998]:    `Incompleteness, Mechanism, and Optimism',
\emph{The Bulletin of Symbolic Logic}, Vol. 4, No. 3, pp. 273-302.

Shapiro, Stewart [2003]:    `Mechanism, truth, and Penrose's new argument',
\emph{Journal of Philosophical Logic}, XXXII(1):19-42.

Solovay, M. R.[1976]: `Provability interpretations of modal logic', \emph{Israel Journal of Mathematics} 25, pp. 287-304.

Stephen, G.Simpson [1985]: \emph{`Harvey Friedman's Research on the Foundations of Mathematics'},
Studies in Logic and the Foundations of Mathematics, vol. 117, North-Holland Publishing,
Amsterdam.

Stephen, G.Simpson (editor) [1987]: `Logic and Combinatorics', in \emph{Contemporary Mathematics}, vol. 65,
AMS, Providence, RI.

Stephen, G.Simpson  [1988]: `Partial Realizations of Hilbert's Program',
\emph{The Journal of Symbolic Logic}, Vol. 53, No. 2, pp. 349-363.

Smith, Peter  [2007]:
\emph{`An Introduction to G\"{o}del's Theorems'}, Cambridge University Press.

Smory\'{n}ski, C. [1977]: `The Incompleteness Theorems', in: J. Barwise (Ed.), Handbook of Mathematical Logic,
North-Holland, Amsterdam,  pp. 821-865.

Stillwell, John [2015]: `What Does ``Depth" Mean in Mathematics'?
\emph{Philosophia Mathematica}, Volume 23, Issue 2,  Pages 215-232.

Tarski, A.,  Mostowski, A.  and Robinson  M. R. [1953]:
\emph{`Undecidable Theories'}, North-Holland.

Urquhart, Alasdair [2015]: `Mathematical Depth',
\emph{Philosophia Mathematica}, Volume 23, Issue 2,  Pages 233-241.

Vaught L. R. [1962]:  `On a theorem of Cobham concerning undecidable theories', In:
Nagel, E., Suppes, P., Tarksi, A. (eds.): Logic, Methodology, and Philosophy of Science,
p. 18. Proceedings of the 1960 International Congress. Stanford, CA: Stanford
University Press.

Visser, Albert [2011]: `Can we make the second incompleteness theorem coordinate
free'? \emph{Journal of Logic and Computation} 21(4), 543-560.

%\bibitem{Visser 14} Albert Visser.  Why the theory $\mathbf{R}$ is special. In Neil Tennant (Eds.), Foundational Adventures-Essays in honour of Harvey M. Friedman (pp. 7-24) (17 p.). College Publication, 2014.

Visser, Albert [2016]: `The Second Incompleteness Theorem:
Reflections and Ruminations', Chapter in G\"{o}del's Disjunction: The scope and limits of mathematical knowledge, edited by Leon Horsten and Philip Welch, Oxford University Press.

%\bibitem{Visser 17} Albert Visser. Another look at the second incompleteness theorem. Logic Group preprint series, 339, (pp. 1-25) (25 p.), 2017.

Visser, Albert [2019]:
`From Tarski to  G\"{o}del: or, how to derive the second incompleteness theorem from the undefinability of truth without self-reference', Reprint.

Vop\v{e}nka, Petr: `A new proof of G\"{o}del's result on non-provability of consistency',
\emph{Bulletin del'Acad\'{e}mie Polonaise des Sciences}.
S\'{e}rie des Sciences Math\'{e}matiques. Astronomiques et Physiques,
14, 111-116.

Willard, E. D.[2001]:  `Self-verifying axiom systems, the incompleteness theorem
and related reflection principles', \emph{Journal of Symbolic Logic}, 66(2):536-596.

Willard, E. D.[2006]: `A generalization of the second incompleteness theorem and
some exceptions to it', \emph{Ann. Pure Appl. Logic}, 141(3):472-496.

%\bibitem{The continuum hypothesis} {Woodin, W. H,} The continuum hypothesis, part I (2001a), Notices of the
%American Mathematical Society 48(6): 567-576.

Zach, Richard [2007]: `Hilbert's Program Then and Now',
Philosophy of Logic, Handbook of the Philosophy of Science, Pages 411-447.

\end{document}